\documentclass[11pt]{amsart}
\usepackage[mathscr]{eucal}
\usepackage{amsfonts}
\usepackage{amsmath}
\usepackage{amsthm}
\usepackage{amssymb}
\usepackage{amscd}
\usepackage{latexsym}
\usepackage{slashbox}
\usepackage{lscape}

\usepackage[pdftex]{graphicx,color}

\theoremstyle{plain}
\newtheorem{definition}[equation]{Definition}

\newtheorem{corollary}[equation]{Corollary}
\newtheorem{lemma}[equation]{Lemma}
\newtheorem{proposition}[equation]{Proposition}
\newtheorem{theorem}[equation]{Theorem}

\theoremstyle{definition}
\newtheorem{remark}[equation]{Remark}

\numberwithin{equation}{subsection}
\newcommand{\Tr}{\text{Tr}}

\title{Hochschild and cyclic (co)homology of preprojective algebras of quivers of type T}
\author{Ching-Hwa Eu}
\address{Department of Mathematics, Massachusetts Institute of
Technology, Cambridge, MA 02139, U.S.A.}
\email{ceu@math.mit.edu}

\begin{document}

\maketitle
\pagestyle{myheadings}
\markboth{Ching-Hwa Eu}{Hochschild and cyclic (co)homology of preprojective algebras of quivers of type T}

\tableofcontents
\section{Introduction}
In this paper, we compute the Hochschild homology spaces and Hochschild cohomology ring of preprojective algebras $A$ of quivers of type T over a field of characteristic zero, together with the grading induced by the natural grading (all arrows have degree $1$). 

This was done for type A (over a field of any characteristic) in \cite{ES2} and for type D/E in \cite{EE2}/\cite{Eu}.

For our computation, we construct a projective resolution of $A$ which turns out to be periodic with period 6.

To compute the cup product, we use the same method as in \cite{ES2} and \cite{Eu}: via the isomorphism $HH^i(A)\equiv \underline{\mathrm{Hom}}(\Omega^iA,A)$ (where for an $A$-bimodule $M$ we write $\Omega M$ for the kernel of its projective cover) we identify elements in $HH^i(A)$ with equivalence classes of maps $\Omega^i(A)\rightarrow A$. For $[f]\in HH^i(A)$ and $[g]\in HH^j(A)$, the product is $[f][g]:=[f\circ\Omega^ig]$ in $HH^{i+j}(A)$. We compute all products $HH^i(A)\times HH^j(A)\rightarrow  HH^{i+j}(A)$ for $0\leq i\leq j\leq 5$. The remaining ones follow from the perodicity of the Schofield resolution and the graded commutativity of the multiplication. 

As in the ADE-case, the preprojective algebra of a quiver of type T is Frobenius. 
But unlike the ADE-case where the Nakayama automorphism has order $2$, in the T-case the Nakayama automorphism is the identity map. Although in the T-case $\Omega^3(A)$ is also the bimodule $A$, where the right action is twisted by an automorphism of order $2$, this is now not the Nakayama automorphism. The consequence is that the projective resolution we get is similar to the Schofield resolution which we used for the ADE quivers, but it involves different A-bimodules.

In Section 6, we recall some properties of $A$. 

We compute the Hochschild homology/cohomology spaces in Section 7: we use the Connes exact sequence and duality theorems from the Calabi-Yau Frobenius property of $A$.

In Section 8, we construct a basis of the Hochschild cohomology space. Then we compute the cup product in Section 9 in terms of this basis.

In Section 10, we use the cup product structure to compute the calculus  ($HH^\bullet(A),HH_\bullet(A)$).

In the Appendix, we include a correction to \cite{ES2} where the cup product for the Hochschild cohomology of preprojective algebras of type $A$ was computed. 

{\bf{Acknowledgements.}} C. Eu wants to thank his advisor P. Etingof and T. Schedler for useful discussions. This work is partially supported by the NSF grant DMS-0504847.

\section{Preliminaries}

\subsection{Quivers and path algebras}
Let $Q$ be a quiver of type $T$ with vertex set $I$ and $|I|=n$. 
We write $x\in Q$ to say that $x$ is an arrow in $Q$. We call the loop $b$. 
Let $Q'=Q\setminus\{b\}$. 

We define $(Q')^*$ to be the quiver obtained from $Q'$ by reversing
all of its arrows. We call $\bar Q'=Q'\cup (Q')^*$ the \emph{double}
of $Q'$. Let $C$ be the adjacency matrix corresponding to the
quiver $\bar Q=\bar Q'\cup\{b\}$. 

The concatenation of arrows generate the \emph{nontrivial
paths} inside the quiver $\bar Q$. We define $e_i$, $i\in I$ to
be the \emph{trivial path} which starts and ends at $i$. The
\emph{path algebra} $P_{\bar Q}=\mathbb{C}\bar Q$ of $\bar Q$ over $\mathbb{C}$ is the $\mathbb{C}$-algebra with basis the paths in $\bar Q$ and the product $xy$ of two paths $x$ and $y$ to be their concatenation if they are compatible and $0$ if not. We define the \emph{Lie bracket} $[x,y]=xy-yx$.

Let $R=\bigoplus_{i\in I}\mathbb{C}e_i$. Then $R$ is a commutative semisimple
algebra, and $P_{\bar Q}$ is naturally an
$R$-bimodule.

\subsection{Frobenius algebras}
Let $A$ be a finite dimensional unital $\mathbb{C}-$algebra. We
call it Frobenius if there is a linear function $f:A\rightarrow\mathbb{C}$, such that the form $(x,y):=f(xy)$ is nondegenerate, or, equivalently, if there exists an isomorphism $\phi:A\stackrel{\simeq}{\rightarrow}A^*$ of left $A-$modules: given $f$, we can define $\phi(a)(b)=f(ba)$, and given $\phi$, we define $f=\phi(1)$.

If $\tilde f$ is another linear function satisfying the same
properties as $f$ from above, then $\tilde f(x)=f(xa)$ for some
invertible $a\in A$. Indeed, we define the form $\{a,b\}=\tilde f(ab)$. Then $\{-,1\}\in A^*$, so there is an $a\in A$, such that $\phi(a)=\{-,1\}$. Then $\tilde f(x)=\{x,1\}=\phi(a)(x)=f(xa)$.

\subsection{The Nakayama automorphism}
Given a Frobenius algebra $A$ (with a function $f$ inducing a
bilinear form $(-,-)$ from above), the automorphism
$\eta:A\rightarrow A$ defined by the equation $(x,y)=(y,\eta(x))$
is called the \emph{Nakayama automorphism} (corresponding to
$f$). 

We note that the freedom in choosing $f$ implies 
that $\eta$ is uniquely determined up to an inner automorphism. Indeed,
let $\tilde f(x)=f(xa)$ and define the bilinear form $\{x,y\}=\tilde f(xy)$. Then
\begin{align*}
\{x,y\}&=\tilde f(xy)=f(xya)=(x,ya)=(ya,\eta(x))=f(ya\eta(x)a^{-1}a)\\
&=(y,a\eta(x)a^{-1}).
\end{align*}

\subsection{The preprojective algebra}
Given a T-quiver $Q$, we define the \emph{preprojective
algebra} $\Pi_Q$ to be the quotient of the path algebra $P_{\bar Q}$ by
the relation $\sum\limits_{a\in Q'}[a,a^*]=b^2$. It is known that 
$\Pi_Q$ is a Frobenius algebra (see \cite{MOV}).

From now on, we write $A=\Pi_Q$.

\subsection{Graded spaces and Hilbert series}

Let $W=\bigoplus_{d\geq0}W(d)$ be a $\mathbb Z_+$-graded
vector space, with finite dimensional homogeneous subspaces. 
We denote by $M[d]$ the same space with
grading shifted by $d$. The graded dual space $M^*$ is defined by the
formula $M^*(d)=M(-d)^*$.  

\begin{definition} \textnormal{(The Hilbert series of vector spaces)}\\
We define the \emph{Hilbert series} $h_W(t)$ to be the series
\begin{displaymath}
h_W(t)=\sum\limits_{d=0}^{\infty}\dim W(d)t^d.
\end{displaymath}
\end{definition}

\begin{definition} \textnormal{(The Hilbert series of bimodules)}\\
Let $W=\bigoplus_{d\geq0}W(d)$ be a $\mathbb{Z_+}$-graded bimodule
over the ring $R$, 
so we can write $W=\bigoplus_{i,j\in I} W_{i,j}$. We define the 
\emph{Hilbert series} $H_W(t)$ to be a matrix valued series with the entries 
\begin{displaymath}
H_W(t)_{i,j}=\sum\limits_{d=0}^{\infty}\dim\ W(d)_{i,j}t^d.
\end{displaymath}
\end{definition}

\section{The calculus structure of the preprojective algebra}
We recall the definition of the calculus.
\subsection{Definition of calculus}
\begin{definition}\emph{(Gerstenhaber algebra)} 
A graded vector space $\mathcal{V}^\bullet$ is a Gerstenhaber algebra if it is equipped with a graded commutative and associative product $\wedge$ of degree 0 and a graded Lie bracket $[,]$ of degree $-1$. These operations have to be compatible in the sense of the following Leibniz rule

\begin{equation}
 [\gamma,\gamma_1\wedge\gamma_2]=[\gamma,\gamma_1]\wedge\gamma_2+(-1)^{k_1(k+1)}\gamma_1\wedge[\gamma,\gamma_2],
\end{equation}
where $\gamma\in\mathcal{V}^k$ and $\gamma_1\in\mathcal{V}^{k_1}$.
\end{definition}
We recall from \cite{CST} that
\begin{definition}\emph{(Precalculus)} 
A precalculus is a pair of a Gerstenhaber algebra $(\mathcal{V}^\bullet,\wedge,[,])$ and a graded vector space $\mathcal{W}^\bullet$ together with 
 \begin{itemize}
  \item a module structure $\iota_\bullet:\mathcal{V}^\bullet\otimes\mathcal{W}^{-\bullet}\rightarrow\mathcal{W}^{-\bullet}$ of the graded commutative algebra $\mathcal{V}^\bullet$ on $\mathcal{W}^{-\bullet}$.
  \item an action $\mathcal{L_\bullet}:\mathcal{V}^{\bullet+1}\otimes\mathcal{W}^{-\bullet}\rightarrow\mathcal{W}^{-\bullet}$ of the graded Lie algebra $\mathcal{V}^{\bullet+1}$ on $\mathcal{W}^{-\bullet}$ which are compatible in the sense of the following equations 
 \begin{equation}
  \iota_a\mathcal{L}_b-(-1)^{|a|(|b|+1)}\mathcal{L}_b\iota_a=\iota_{[a,b]},
 \end{equation}
and 
\begin{equation}
 \mathcal{L}_{a\wedge b}=\mathcal{L}_a\iota_b+(-1)^{|a|}\iota_a\mathcal{L}_b.
\end{equation}
 \end{itemize}
\end{definition}
\begin{definition}\emph{(Calculus)} 
 A calculus is a precalculus $(\mathcal{V}^\bullet,\mathcal{W}^\bullet,[,],\wedge,\iota_\bullet,\mathcal{L}_\bullet)$ with a degree 1 differential $d$ on $\mathcal{W}^\bullet$ such that the \emph{Cartan identity},
 \begin{equation}\label{Cartan}
  \mathcal{L}_a=d\iota_a-(-1)^{|a|}\iota_ad,
 \end{equation}
 holds.
\end{definition}
Let $A$ be an associative algebra. The contraction of the Hochschild cochain $P\in C^k(A,A)$ with the Hochschild chain $(a_0,a_1,\ldots,a_n)$ is defined by
 
\begin{equation}
 I_P(a_0,a_1,\ldots,a_n)=
 \left\{
 \begin{array}{cc}
(a_0P(a_1,\ldots,a_k),a_{k+1},\ldots,a_n) & n\geq k,\\
0   &   \mbox{else}.
\end{array}
\right.
\end{equation}

We have
\begin{proposition}\emph{\textbf{(Yu. Daletski, I. Gelfand and B. Tsygan \cite{DGT})}}
 The contraction $I_P$ together with the Connes differential, the Gerstenhaber bracket, the cup product and the action of cochains on chains (\cite[(3.5), page 46]{D}) induce on the pair $(HH^\bullet(A,A), HH_\bullet(A,A))$ a structure of calculus.
\end{proposition}

\section{Results about Hochschild and cyclic (co)homology of $A$}
\begin{definition}
 We define the spaces
 \begin{align*}
  U&=\bigoplus_{d<h-2}HH^0(A)(d)[2]\quad and\\
  K&=HH^2(A)[2].
 \end{align*}
\end{definition}
\begin{theorem}\label{UK}
The spaces $U$ and $K$ have the following properties:
\begin{description}
 \item[(a)] $U$ has Hilbert series
 \begin{equation}
  h_{U}(t)=\sum\limits_{i=0}^{n-1}t^{2i}.
 \end{equation}
 \item[(b)] $K$ is n-dimensional and sits in degree zero.
\end{description}
\end{theorem}

\begin{theorem}[Hochschild cohomology]\label{cohomology}
 The Hochschild cohomology spaces are given by 
 \begin{eqnarray*}
 HH^0(A)&=&U[-2]\oplus R^*[h-2],\\
 HH^1(A)&=&U[-2],\\
 HH^2(A)&=&K[-2],\\
 HH^3(A)&=&K^*[-2],\\
 HH^4(A)&=&U^*[-2],\\
 HH^5(A)&=&U^*[-2],\\
 HH^6(A)&=&U[-2h-2],\\
 HH^{6k+i}(A)&=&HH^i(A)[-2kh]\quad\forall{i\geq 1}.
\end{eqnarray*}
\end{theorem}

\begin{theorem}[Hochschild homology]\label{homology}
 The Hochschild homology spaces are given by
\begin{eqnarray*}
 HH_0(A)&=&U^*[h]\oplus R,\\
 HH_1(A)&=&U^*[h],\\
 HH_2(A)&=&K^*[h],\\
 HH_3(A)&=&K[h],\\
 HH_4(A)&=&U[h],\\
 HH_5(A)&=&U[h],\\
 HH_6(A)&=&U^*[3h],\\
 HH_{6k+i}(A)&=&HH_i(A)[2kh]\quad\forall i\geq1.
\end{eqnarray*}
\end{theorem}

\begin{theorem}[Cyclic homology]\label{cyclic}
 The cyclic homology spaces are given by
\begin{eqnarray*}
 HC_0(A)&=&U^*[h]\oplus R,\\
 HC_1(A)&=&0,\\
 HC_2(A)&=&K^*[h],\\
 HC_3(A)&=&0,\\
 HC_4(A)&=&U[h],\\
 HC_5(A)&=&0,\\
 HC_6(A)&=&U^*[3h],\\
 HC_{6k+i}(A)&=&HH_i(A)[2kh]\quad\forall i\geq1.
\end{eqnarray*}
\end{theorem}

Let $(U[-2])_+$ be the positive degree part of $U[-2]$ (which lies in non-negative degrees). 

We have a decomposition $HH^0(A)=\mathbb{C}\oplus (U[-2])_+\oplus L[-h-2]$ where we have the natural identification $(U[-2])(0)=\mathbb{C}$. This identification also gives us a decomposition $HH^*(A)=\mathbb{C}\oplus HH^*(A)_+$.\\
We also decompose $U=U^{top}\oplus U_-$, where $U^{top}$ is the top degree part of $U$ and a one-dimesional space.

We give a brief description of the product structure in $HH^*(A)$ which will be computed in this paper. Since the product $HH^i(A)\times HH^j(A)\rightarrow HH^{i+j}(A)$ is graded-commutative, we can assume $i\leq j$ here.

Let $z_0=1\in\mathbb{C}\subset U[-2]\subset HH^0(A)$ (in lowest degree $0$), \\
$\theta_0$ the corresponding element in $HH^1(A)$ (in lowest degree $0$),\\
$\psi_0$ the dual element of $z_0$ in $U^*[-2]\subset HH^5(A)$ (in highest degree $-4$), i.e. $\psi_0(z_0)=1$, \\
$\zeta_0$ the corresponding element in $U^*[-2]\subset HH^4(A)$ (in highest degree $-4$), that is the dual element of $\theta_0$, $\zeta_0(\theta_0)=1$,\\
$\varphi_0:HH^0(A)\rightarrow HH^6(A)$
the natural quotient map (which induces the natural isomorphism $U[-2]\rightarrow U[-2h-2])$.\\

\begin{theorem}[Cup product]\label{cup}
 \begin{enumerate}
\item
The multiplication by $\varphi_0(z_0)$ induces the natural isomorphisms\\ $\varphi_i:HH^i(A)\rightarrow HH^{i+6}(A)$ $\forall i\geq 1$ and the natural quotient map $\varphi_0$. Therefore, it is enough to compute products $HH^i(A)\times HH^j(A)\rightarrow HH^{i+j}(A)$ with $0\leq i\leq j\leq 5$.
\item The $HH^0(A)$-action on $HH^i(A)$.
\begin{enumerate}
\item ($(U[-2])_+$-action).\\ The action of $(U[-2])_+$ on $U[-2]\subset HH^1(A)$ corresponds to the multiplication 
\begin{eqnarray*}
(U[-2])_+\times U[-2]&\rightarrow& U[-2], \\
(u,v)&\mapsto& u\cdot v
\end{eqnarray*}
in  $HH^0(A)$, projected on $U[-2]\subset HH^0(A)$.\\
$(U[-2])_+$ acts on $U^*[-2]=HH^4(A)$ and $U^*[-2]\subset HH^5(A)$ the following way:
\begin{eqnarray*}
(U[-2])_+\times U^*[-2]&\rightarrow& U^*[-2],\\
(u,f)&\mapsto&u\circ f,
\end{eqnarray*}
where $(u\circ f)(v)=f(uv)$.\\
$(U[-2])_+$ acts by zero on $R^*[h-2]\subset HH^0(A)$, $HH^2(A)$ and $HH^3(A)$.
\item ($R^*[h-2]$-action).\\
$R^*[h-2]$ acts by zero on $HH^*(A)_+$.
\end{enumerate}
\item (Zero products).\\
For all odd $i,j$, the cup product $HH^i(A)\cup HH^j(A)$ is zero.
\item ($HH^1(A)$-products).
\begin{enumerate}
\item The multiplication \[HH^1(A)\times HH^4(A)=U[-2]\times U^*[-2]\rightarrow HH^5(A)\] is the same one as the restriction of \[HH^0(A)\times HH^5(A)\rightarrow HH^5(A)\] on $U[-2]\times U^*[-2]$.\\
\item The multiplication of the subspace $U[-2]_+\subset HH^1(A)$ with $HH^2(A)$ is zero.\\
\item The multiplication by $\theta_0$ induces a symmetric isomorphism \[\alpha:HH^2(A)=K[-2]\rightarrow K^*[-2]=HH^3(A),\] given by the matrix $(2n+1)(2-C')^{-1}$, where $C'$ is obtained from the adjacency matrix by changing the sign on the diagonal.
\end{enumerate}
\item ($HH^2(A)$-products).
\begin{eqnarray*}
HH^2(A)\times HH^2(A)&\rightarrow& HH^4(A),\\ 
(a,b)&\mapsto& \langle-,-\rangle\zeta_0
\end{eqnarray*} 
is given by $\langle-,-\rangle=\alpha$ where $\alpha$ is regarded as a symmetric bilinear form.\\

$HH^2(A)\times HH^3(A)\rightarrow HH^5(A)$ is the multiplication
\begin{eqnarray*}
K[-2]\times K^*[-2]&\rightarrow& HH^5(A),\\
(a,y)&\mapsto&y(a)\psi_0.
\end{eqnarray*}
\item (Products involving $U^*[-2]$).
\begin{enumerate}
\item (($(U_-)^*[-2]$-action).\\
$(U_-)^*[-2]\subset HH^i(A)$, $i=4,5$ acts by zero on $HH^j(A)$, $j=2,3,4,5$.
\item We can choose a $\zeta'\in(U^{top})^*[-2]\in HH^4(A)$, such that for the corresponding elements $z'\in U^{top}[-2]\subset HH^0(A)$, $\theta'\in U^{top}[-2]\subset HH^1(A)$ and $\psi'\in (U^{top})^*[-2]\subset HH^5(A)$ we get the formulas below.
\begin{enumerate}
\item $HH^4(A)\times HH^2(A)\rightarrow HH^6(A)$. The multiplication with $\zeta'$ gives us a map
\begin{eqnarray*}
 K[-2]=\mathbb{C}^I&\rightarrow&U^{top}[-2h-2],\\
 v=(v_i)_{i\in I}&\mapsto& v_s \varphi_0(z').
\end{eqnarray*}
\item $HH^5(A)\times HH^2(A)\rightarrow HH^7(A)$. This pairing
\[
 U^*[-2]\times K[-2]\rightarrow U[-2h-2]
\]
is the same as the corresponding pairing 
\[
 HH^4(A)\times HH^2(A)\rightarrow HH^6(A).
\]
\item $HH^4(A)\times HH^3(A)\rightarrow HH^7(A)$. The multiplication with $\zeta'$ gives us a map
\begin{eqnarray*}
 K^*[-2]=\mathbb{C}^I&\rightarrow&U^{top}[-2h-2],\\
 w=(w_i)_{i\in I}&\mapsto&\sum\limits_{i\in I}(n-d(i,s)) w_i \varphi_0(\theta'). 
\end{eqnarray*}
where $d(i,j)$ is the distance between two vertices $i,\,j$.
\item $HH^4(A)\times HH^4(A)\rightarrow HH^8(A)$. $\zeta'^2$ gives us the vector $(\delta_{is})_{i\in I}$ in $K[-2h-2]=HH^8(A)$.
\item $HH^4(A)\times HH^5(A)\rightarrow HH^9(A)$. $\zeta'^2\psi'^2$ gives us the vector $(n-d(i,s))_{i\in I}$ in $K^*[-2h-2]=HH^9(A)$. $HH^4(A)$ annihilates $(U_-)^*[-2]\subset HH^5(A)$.
\end{enumerate}
\end{enumerate}
\end{enumerate}
\end{theorem}

Comparing this theorem with the results about the $A$-case in \cite{ES2}, we get the following:

\begin{corollary}[relation to the A-case]\label{A}
 Let $\omega_1,\ldots,\omega_n$ be a basis of $R^*[h-2]\subset HH^0(A)$. Then we have
 \begin{equation}
  HH^*(\Pi_{T_n})=HH^*(\Pi_{A_{2n}})[\omega_1,\ldots,\omega_n]/(R^*[h-2] HH^*(\Pi_{A_{2n}})_+).
 \end{equation}
We can write $HH^*(\Pi_{A_{2n}})$ as a quotient
\begin{equation}
HH^*(\Pi_{A_{2n}})=HH^*(\Pi_{T_n})/(R^*[h-2]).
\end{equation}
\end{corollary}

\section{Results about the Calculus}
We will introduce for every $m\geq0$ an isomorphism
\begin{equation}
\mathbb{D}:HH_m(A)\stackrel{\sim}{\rightarrow}HH^{6m+5}(A)[(2m+1)h+2]
\end{equation}
 which intertwines contraction and cup-product maps. 

In Section \ref{basis}, we will introduce basis elements $z_k\in U[-2]\subset HH^0(A)$, $\theta_k\in HH^1(A)$, $f_k\in HH^2(A)$, $h_k\in HH^3(A)$, $\zeta_k\in HH^4(A)$ and $\psi_k\in HH^5(A)$.

For $c_k\in HH^i(A)$, $0\leq i\leq 5$, we write $c_k^{(s)}$ for the corresponding cocycle in $HH^{i+6s}$. 
We write $c_{k,t}$ for a cycle in $HH_{j+6t}$, $1\leq j\leq6$ which equals $\mathbb{D}^{-1}(c_k^{(s)})$.

The map $\alpha: K\rightarrow K^*$, given by a matrix $M_\alpha$, is introduced in Subsection \ref{H1xH2}.

We state the results in terms of these bases of $HH^\bullet(A)$ and $HH_\bullet(A)$.

\begin{theorem}
 The calculus structure is given by tables 1, 2, 3 and the \emph{Connes differential} B, given by
\begin{eqnarray*}
B_{6s}(\psi_{k,s})&=&((2s+1)h-2-k)\zeta_{k,s},\\
B_{1+6s}&=&0,\\
B_{2+6s}(h_{k,s})&=&(2s+1)h\alpha^{-1}(h_{k,s}),\\
B_{3+6s}&=&0,\\
B_{4+6s}(\theta_{k,s})&=&((2s+1)h+2+k)z_{k,s},\\
B_{5+6s}&=&0.
\end{eqnarray*}
\end{theorem}

\begin{landscape}
\thispagestyle{empty}
\begin{table}
\begin{tabular}[b]{|c||c|c|c|c|c|c|c|c|}
\hline
\backslashbox{$a$}{$b$}&$\psi_{l,t}$&$\zeta_{l,t}$&$h_{l,t}$&$f_{l,t}$&$\theta_{l,t}$&$z_{l,t}$\\\hline\hline
$z_k^{(s)}$&$(z_k\psi_{l})_{t-s}$&$(z_k\zeta_{l})_{t-s}$&$\delta_{k0}h_{l,t-s}$&$\delta_{k0}f_{l,t-s}$&$(z_k\theta_l)_{t-s}$&$(z_kz_l)_{t-s}$\\\hline
$\omega_k$&$0$&$0$&$0$&$0$&$0$&$0$\\\hline
$\theta_k^{(s)}$&$0$&$(z_k\psi_l)_{t-s}$&$0$&$\delta_{k0}\alpha(f_{l,t-s})$&$0$&$(z_l\theta_k)_{t-s}$\\\hline
$f_k^{(s)}$&$\begin{array}{ll}\delta_{l,h-3}k\cdot\\\theta_{h-3,t-s-1}\end{array}$&$\begin{array}{ll}\delta_{l,h-3}k\cdot\\z_{l,t-s-1}\end{array}$&$\delta_{kl}\psi_{0,t-s}$&$(M_\alpha)_{kl}\zeta_{0,t-s}$&$\delta_{l0}\alpha(f_{k,t-s})$&$\delta_{l0}f_{k,t-s}$\\\hline
$h_k^{(s)}$&$0$&$\begin{array}{ll}\delta_{k,n}\delta_{l,h-3}\cdot\\\theta_{h-3,t-s-1}\end{array}$&$0$&$\delta_{kl}\psi_{0,t-s}$&$0$&$\delta_{l0}h_{k,t-s}$\\\hline
$\zeta_k^{(s)}$&$\begin{array}{ll}\delta_{k,h-3}\delta_{l,h-3}\cdot\\\alpha(f_{n,t-s-1})\end{array}$&$\begin{array}{ll}\delta_{k,h-3}\delta_{l,h-3}\cdot\\ f_{n,t-s-1}\end{array}$&$\begin{array}{ll}\delta_{k,h-3}\delta_{l,n}\cdot\\\theta_{k,t-s-1}\end{array}$&$\begin{array}{ll}\delta_{k,h-3}l\cdot\\z_{k,t-s-1}\end{array}$&$(z_l\psi_k)_{t-s}$&$(z_l\zeta_k)_{t-s}$\\\hline
$\psi_k^{(s)}$&$0$&$\begin{array}{ll}\delta_{k,h-3}\delta_{l,h-3}\cdot\\\alpha(f_{n,t-s-1})\end{array}$&$0$&$\begin{array}{ll}\delta_{k,h-3}l\cdot\\\theta_{h-3,t-s-1}\end{array}$&$0$&$(z_k\psi_l)_{t-s}$\\\hline
\end{tabular}
\caption{contraction map $\iota_a(b)$}
\label{contraction}
\end{table}

\end{landscape}

\begin{landscape}
\begin{table}
\begin{center}
\begin{tabular}[b]{|c||@{}c@{\hspace{0 cm}}|c|@{\hspace{-0.18 cm}}c@{\hspace{-0.18 cm}}|@{\hspace{-0.18 cm}}c@{\hspace{-0.18 cm}}|@{\hspace{-0.1 cm}}c@{\hspace{-0.18 cm}}|@{\hspace{-0.18 cm}}c@{\hspace{-0.18 cm}}|@{\hspace{-0.1 cm}}c@{\hspace{-0.18 cm}}|}
\hline
\backslashbox{$a$}{$b$}&$z_l^{(t)}$&$\omega_l$&$\theta_l^{(t)}$&$f_l^{(t)}$&$h_l^{(t)}$&$\zeta_l^{(t)}$&$\psi_l^{(t)}$\\\hline\hline
$z_k^{(s)}$&$0$&$0$&$\begin{array}{cc}(k-2sh)\cdot\\(z_kz_l)^{(s+t)}\end{array}$&$0$&$\begin{array}{cc}-2\delta_{k0}sh\cdot\\\alpha^{-1}(h_l^{(s+t)})\end{array}$&$0$&$\begin{array}{cc}(k-2sh)\cdot\\(z_k\zeta_{l})^{(s+t)}\end{array}$\\\hline
$\omega_k$&&$0$&$0$&$0$&$0$&$0$&$0$\\\hline                                                                      $\theta_k^{(s)}$&&&$\begin{array}{cc}(l-k+2(s-t)h)\cdot\\(z_k\theta_{l})^{(s+t)}\end{array}$&$\begin{array}{cc}-2(1+th)\cdot\\\delta_{k0}f_l^{(s+t)}\end{array}$&$\begin{array}{cc}2(-1+(s-t)h)\cdot\\\delta_{k0}h_l^{(s+t)}\end{array}$&$\begin{array}{cc}-(4+l+2th)\cdot\\(z_k\zeta_{l})^{(s+t)}\end{array}$&$\begin{array}{cc}-(4+k+l+2(t-s)h)\cdot\\(z_k\psi_{l})^{(s+t)}\end{array}$\\\hline
$f_k^{(s)}$&&&&$0$&$\begin{array}{cc}-2(1+sh)\cdot\\\delta_{kl}\zeta_0^{(s+t)}\end{array}$&$0$&$\begin{array}{cc}-2(k+1)\cdot\\ (1+sh)\cdot\\\delta_{l,h-3}z_{h-3}^{(s+t+1)}\end{array}$\\\hline
$h_k^{(s)}$&&&&&$\begin{array}{c}2(s-t)h\cdot\\(M_\alpha^{-1})_{kl}\cdot\\\psi_0^{(s+t)}\end{array}$&$\begin{array}{cc}-(h+1+2th)\cdot\\\delta_{k,\frac{h-3}{2}}\delta_{l,h-3}\cdot\\z_{h-3}^{(s+t+1)}\end{array}$&$\begin{array}{c}(2(s-t)h-(h-1))\cdot\\\delta_{k,\frac{h-3}{2}}\delta_{l,h-3}\cdot\\\theta_{h-3}^{(s+t+1)}\end{array}$\\\hline
$\zeta_k^{(s)}$&&&&&&$0$&$\begin{array}{c}-(2sh+h+1)\cdot\\\delta_{k,h-3}\delta_{l,h-3}\cdot\\ f_{\frac{h-3}{2}}^{(s+t+1)}\end{array}$\\\hline
$\psi_k^{(s)}$&&&&&&&$\begin{array}{c}2(s-t)h\cdot\\\delta_{k,h-3}\delta_{l,h-3}\cdot\\\alpha(f_{\frac{h-3}{2}}^{(s+t+1)})\end{array}$\\\hline
\end{tabular}
\end{center}
\caption{Gerstenhaber bracket $[a,b]$}
\label{Gerstenhaber bracket}
\end{table}
\end{landscape}

\begin{landscape}
\thispagestyle{empty}	
 \begin{table}
 \begin{tabular}[b]{|@{\hspace{0.1 cm}}c@{\hspace{0.1 cm}}||@{\hspace{-0.18 cm}}c@{\hspace{-0.18 cm}}|@{\hspace{-0.18 cm}}c@{\hspace{-0.18 cm}}|@{\hspace{-0.18 cm}}c@{\hspace{-0.18 cm}}|@{\hspace{-0.18 cm}}c@{\hspace{-0.18 cm}}|@{\hspace{-0.18 cm}}c@{\hspace{-0.18 cm}}|@{\hspace{-0.18 cm}}c@{\hspace{-0.18 cm}}|}
\hline
\backslashbox{$a$}{$b$}&$\psi_{l,t}$&$\zeta_{l,t}$&$h_{l,t}$&$f_{l,t}$&$\theta_{l,t}$&$z_{l,t}$\\\hline\hline
$\theta_k^{(s)}$&$\begin{array}{c}((2t+1)h-2-l)\cdot\\(z_k\psi_l)_{t-s}\end{array}$&$\begin{array}{c}((2(t-s)+1)h\\-2-l+k)\\(z_k\zeta_l)_{t-s}\end{array}$&$\begin{array}{c}(2t+1)h\cdot\\\delta_{k0}h_{l,t-s}\end{array}$&$\begin{array}{c}(2(t-s)+1)h\cdot\\\delta_{k0}f_{l,t-s}\end{array}$&$\begin{array}{c}((2t+1)h+2+l)\cdot\\(z_k\theta_l)_{t-s}\end{array}$&$\begin{array}{c}((2(t-s)+1)h+2+k+l)\cdot\\(z_kz_l)_{t-s}\end{array}$\\\hline
$f_k^{(s)}$&$\begin{array}{c}-2k(1+sh)\cdot\\\delta_{l,h-3}z_{h-3,t-s-1}\end{array}$&$0$&$\begin{array}{c}-2(1+sh)\cdot\\\delta_{kl}\zeta_{0,t-s}\end{array}$&$0$&$\begin{array}{c}-2(1+sh)\cdot\\\delta_{l0}f_{k,t-s}\end{array}$&$0$\\\hline
$h_k^{(s)}$&$\begin{array}{c}(2th+1)\cdot\\\delta_{k,n}\delta_{l,h-3}\cdot\\\theta_{h-3,t-s-1}\end{array}$&$\begin{array}{c}(2(t-s)h-1)\cdot\\\delta_{k,n}\delta_{l,h-3}\cdot\\z_{h-3,t-s-1}\end{array}$&$\begin{array}{c}(2t+1)h\cdot\\(M_\alpha^{-1})_{lk}\cdot\\\psi_{0,t-s}\end{array}$&$\begin{array}{c}(2(t-s+1)h-2)\cdot\\\delta_{kl}\zeta_{0,t-s}\end{array}$&$\begin{array}{c}((2t+1)h+2)\cdot\\\delta_{l0}h_{k,t-s}\end{array}$&$\begin{array}{c}\delta_{l0}(2(t-s)+1)h\cdot\\\alpha^{-1}(h_{k,t-s})\end{array}$\\\hline
$\zeta_k^{(s)}$&$\begin{array}{c}-((2s+1)h+1)\cdot\\\delta_{k,h-3}\delta_{l,h-3}\cdot\\f_{n,t-s-1}\end{array}$&$0$&$\begin{array}{c}-((2s+1)h+1)\cdot\\\delta_{k,h-3}\delta_{l,n}\cdot\\z_{h-3,t-s-1}\end{array}$&$0$&$\begin{array}{c}-(2sh+4+k)\cdot\\(z_l\zeta_k)_{t-s}\end{array}$&$0$\\\hline
$\psi_k^{(s)}$&$\begin{array}{c}(2th+1)\cdot\\\delta_{k,h-3}\delta_{l,h-3}\cdot\\\alpha(f_{n,t-s-1})\end{array}$&$\begin{array}{c}(2(t-s)-1)h\cdot\\\delta_{k,h-3}\delta_{l,h-3}\cdot\\ f_{n,t-s-1}\end{array}$&$\begin{array}{c}(2t+1)h\cdot\\\delta_{k,h-3}\delta_{l,n}\cdot\\\theta_{h-3,t-s-1}\end{array}$&$\begin{array}{c}l\cdot\\
(2(t-s)h\\-1)\cdot\\\delta_{k,h-3}\\z_{h-3,t-s-1}\end{array}$&$\begin{array}{c}((2t+1)h\\+2+l)\cdot\\(z_l\psi_k)_{t-s}\end{array}$&$\begin{array}{c}((2(t-s)+1)h\\-2-k+l)\cdot\\(z_l\zeta_k)_{t-s}\end{array}$\\\hline
$z_k^{(s)}$&$\begin{array}{c}(k-2sh)\cdot\\(z_k\zeta_l)_{t-s}\end{array}$&$0$&$\begin{array}{c}(k-2sh)\cdot\\\alpha^{-1}(h_{l,t-s})\end{array}$&$0$&$\begin{array}{c}(k-2sh)\\(z_k\theta_l)_{t-s}\end{array}$&$0$\\\hline
\end{tabular}
\caption{Lie derivative $\mathcal{L}_a(b)$}
\label{Lie derivative}
\end{table}
\end{landscape}

\section{Properties of $A$}

\subsection{Labeling}
We choose a labeling of the quiver $T_n$.

\begin{figure}[htp]
\centerline{
\includegraphics[angle=-90]{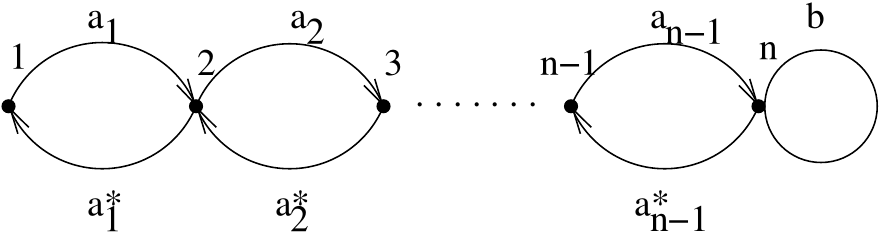}}
\caption{$T_n$ quiver}
\end{figure}

\subsection{Bases and Hilbert series}\label{bases}
From \cite{MOV}, we know that $H_A(t)=(1+t^h)(1-Ct+t^2)^{-1}$, where $C$ is the adjacency matrix of the underlying graph. By choosing the labeling of the quiver above, we get

\begin{equation}
 (\dim e_iAe_j)_{i,j\in I}=H_A(1)=2\cdot\left[
 \begin{array}{ccccccc}
  1&1&1&\cdots&\cdots& 1&1\\
  1&2&2&\cdots&\cdots& 2&2\\
  1&2&3&\cdots&\cdots& 3&3\\
  \vdots&\vdots&\vdots&\ddots&&\vdots&\vdots\\
  \vdots&\vdots&\vdots&&\ddots&\vdots&\vdots\\
  1&2&3&\cdots&\cdots&n-1&n-1\\
  1&2&3&\cdots&\cdots&n-1&n
 \end{array}
 \right].
\end{equation}

We will work with explicit bases $\mathfrak{B}_i$ of $e_iAe_i$. 
The $i^{th}$ diagonal entry of $H_A(t)$ is $\sum\limits_{j=0}^{i-1}t^{2j}+t^{2n-1-2j}$, and since in $A$ all paths starting and ending in the same vertex with the same length are equivalent, we can say that bases of $\mathfrak{B}_i$ are given by paths of length $2j,\,0\leq j\leq i-1$ and $2n-1-2j,\,0\leq j\leq i-1$ (one of each length). We call $c_{i,k}$ to be a nonzero element in $A$, represented by a path of length $k$, starting and ending at $i$.

\subsection{The trace function}
For the $T$-quiver, the Nakayama automorphism is trivial. The bilinear form $(-,-)$ which comes with our Frobenius algebra $A$ is given by a trace function $\Tr:A\rightarrow\mathbb{C}$ of degree $-(2n-1)$ by\\ $(x,y)=\Tr(xy)$. We work with an explicit trace function which maps a polynomial of degree $2n-1$ to the sum of its coefficients. Given the basis $(e_i)_{i\in I}$, we denote $(\omega_i)_{i\in I}$ as its dual basis, i.e. $\omega_i$ is the monomial of top degree in $e_iAe_i$ with coefficient $1$.

\subsection{The quotient $A/[A,A]$}
The quotient $A/[A,A]$ turns out to be different than in the ADE-case.

\begin{proposition}\label{A/[A,A]}
 The quotient is 
\[
 A/[A,A]=R\oplus\langle b^i|\,i\text{odd}\rangle.
\]
\end{proposition}
\begin{proof}
The commutator $[A,A]$ is the linear span of 
\begin{itemize}
\item paths $p_{ki}=[p_{ki},e_i]$ from $i$ to $k$, $i\neq k$ and 
\item $p_{ii}-p_{jj}=[p_{ij},p_{ji}]$, where $p_{kl} \in e_kAe_l$, i.e. all differences $p_{ii}-p_{jj}$ where $p_{ii}$ and $p_{jj}$ have same degree $>0$. Since all paths in $e_1Ae_1$ of degree $>0$ give us a zero element in $A$, this gives us that \\$e_iAe_i(d)\subset[A,A]$ $\forall i\in I$ and even $d>0$.
\end{itemize}
From above, we get $e_iAe_k\subset [A,A]$. 
Since all paths in $e_1Ae_1$ of even degree $>0$ give us a zero element in $A$ and $p_{ii}-p_{jj}\in[A,A]$ for any pair of paths of same degree, this implies that $e_iAe_i(d)\subset[A,A]$ $\forall i\in I$ and for all even $d>0$.

$R\cap[A,A]=0$ since $R$ is a commutative ring.

So the quotient $A/[A,A]$ is spanned only by $R$ and by odd degree paths $p_{ii}\in e_iAe_i$, and the only relations involving those is $p_{ii}=p_{jj}$. Let $p_{ii}$ have odd degree $d\leq 2n-1$, then $(a_i^*a_i)^{\frac{2n-1-d}{2}}p_{ii}=\omega_i$, so $p_{ii}\neq 0$ in $A$. So if take one path $p_{ii}\in e_iAe_i$ (for some $i\in I$) in each odd degree $\leq 2n-1$, we get a basis in $A/([A,A]+R)$. Specifically we can choose odd powers of $b$ as a basis.
\end{proof}

\section{Hochschild and cyclic (co)homology of $A$}
In this section, we prove Theorems \ref{cohomology} and \ref{homology}: we construct a projective resolution of $A$, prove duality theorems and compute the Hochschild and cyclic cohomology/homology spaces.
\subsection{A periodic projective resolution of $A$}
Let $\phi$ be an automorphism of $A$, such that $\phi(a)=a\,\forall a\in\bar Q'$ and $\phi(b)=-b$. Note that $\phi=-$Id on $A^{top}$. 
Define the 
$A-$bimodule ${}_1A_\phi$ obtained from $A$ by twisting the
right action by $\phi$, i.e. 
${}_1A_\phi=A$ as a vector space, 
and $\forall x,z\in A,y\in {}_1A_\phi:x\cdot y\cdot z=xy\phi(z).$
Introduce the notation $\epsilon_a=1$ if $a\in Q'$,
$\epsilon_a=-1$ if $a\in (Q')^*$, and let $\epsilon_b=1$. Let $x_i$ be a homogeneous basis $\mathfrak{B}$ 
of $A$ and $x_i^*$
the dual basis under the form attached to the Frobenius algebra
$A$. Let $V$ be the bimodule spanned by the edges of $\bar Q$. 

We start with the following complex:
\[
S_\bullet:\quad 0\rightarrow {}_1A_\phi[h]\stackrel{i}{\rightarrow}A\otimes_R A[2]\stackrel{d_2}{\rightarrow}A\otimes_RV\otimes_RA\stackrel{d_1}{\rightarrow}A\otimes_RA\stackrel{d_0}{\rightarrow}A\rightarrow 0,
\]
where
\begin{align*}
d_0(x\otimes y)&=xy,\\
d_1(x\otimes v\otimes y)&=xv\otimes y-x\otimes vy,\\
d_2(z\otimes t)&=\sum\limits_{i=1}^n\epsilon_{a_i}za_i\otimes a_i^*\otimes t+\sum\limits_{i=1}^n\epsilon_{a_i}z\otimes a_i\otimes a_i^*t\\
&\quad-zb\otimes b\otimes t-z\otimes b\otimes bt,\\
i(x)&=x\sum\limits_{x_i\in\mathfrak{B}} \phi(x_i)\otimes x_i^*.
\end{align*}

$d_id_{i+1}=0$ for $i=0,1,2$ is obvious. We show $d_2i=0$:
We have
\begin{eqnarray*}
 d_2(i(1))&=&d_2(\sum\limits_{x_i\in \mathfrak{B}}\phi(x_i)\otimes x_i^*)\\
 &=&\sum\limits_{x_i\in \mathfrak{B}}\sum\limits_{j=1}^n\epsilon_{a_j}\phi(x_i)a_j\otimes a_j^*\otimes x_i^*+\sum\limits_{x_i\in \mathfrak{B}}\sum\limits_{j=1}^n\epsilon_{a_j}\phi(x_i)\otimes a_j\otimes a_j^*x_i^*\\
 &&-\sum\limits_{x_i\in \mathfrak{B}}\phi(x_i)b\otimes b\otimes x_i^*-\sum\limits_{x_i\in\mathfrak{ B}}\phi(x_i)\otimes b\otimes bx_i^*.
\end{eqnarray*}
The first two terms cancel since 
\begin{eqnarray*}
 \forall a\in\bar Q':\quad\sum\limits_{x_i\in \mathfrak{B}}\phi(x_i)a\otimes a^*\otimes x_i^*&=&\sum\limits_{x_i,x_j\in \mathfrak{B}}(\phi(x_i)a,-\phi(x_j^*))\phi(x_j)\otimes a^*\otimes x_i^*\\&=&\sum\limits_{x_i,x_j\in \mathfrak{B}}\phi(x_j)\otimes a^*\otimes(ax_j^*,x_i) x_i^*\\&=&\sum\limits_{x_i\in \mathfrak{B}}\phi(x_i)\otimes a^*\otimes ax_i^*.
\end{eqnarray*}
The last two terms cancel since 

\begin{eqnarray*}
\sum\limits_{x_i\in \mathfrak{B}}\phi(x_i)b\otimes b\otimes x_i^*&=&
\sum\limits_{x_i,x_j\in \mathfrak{B}}(\phi(x_i)b, -\phi(x_j^*))\phi(x_j)\otimes b\otimes x_i^*\\
&=&\sum\limits_{x_i,x_j\in\mathfrak{B}}\phi(x_j)\otimes b\otimes (-bx_j^*,x_i)x_i^*\\
&=&-\sum\limits_{x_i\in\mathfrak{B}}\phi(x_i)\otimes b\otimes bx_i^*.
\end{eqnarray*}

\begin{lemma}
 $S_\bullet$ is self dual.
\end{lemma}
\begin{proof}
 We introduce the nondegenerate forms
 \begin{itemize}
  \item $(x,y)_\phi=\Tr(x\phi(y))$ on $A$,
  \item $(x\otimes x',y\otimes y')_\phi=\Tr(x\phi(y'))\Tr(x'y)$ on $A\otimes_R A$ and
  \item $(x\otimes\alpha\otimes x',y\otimes\beta\otimes y')_\phi=\Tr(x\phi(y'))\Tr(x'y)(\alpha,\beta)$ on 
  $A\otimes_RV\otimes_RA$, where we define the form on $V$ to be $(\alpha,\beta)=\delta_{\alpha^*,\beta}\epsilon_\beta$.
 \end{itemize}
We apply the functor $(-)^\star=\text{Hom}_\mathbb{C}(-,\mathbb{C})$ and make the identifications\\ $A^\star\simeq A$, $(A\otimes_RA)^\star\simeq A\otimes_RA$ and $(A\otimes_RV\otimes_RA)^\star\simeq A\otimes_RV\otimes_RA$ by the map $x\mapsto(-,x)_\phi$.

We have
\begin{align*}
 (i(x),y\otimes z)_\phi&=(\sum\limits_{x_i\in\mathfrak{B}}x\phi(x_i)\otimes x_i^*,y\otimes z)_\phi=\sum\limits_{x_i\in\mathfrak{B}}\Tr(x\phi(x_i)\phi(z))\Tr(x_i^*y)\\
 &=-\sum\limits_{x_i\in\mathfrak{B}}\Tr(\phi(x)x_iz)\Tr(x_i^*y)=-\sum\limits_{x_i\in\mathfrak{B}}\Tr(z\phi(x)x_i)\Tr(x_i^*y)\\
 &=-\Tr(z\phi(x)y)=\Tr(x\phi(yz))=(x,yz)_\phi\\
 &=(x,d_0(y\otimes z)),
\end{align*}
so $i=d_0^\star$.

We have
\begin{align*}
 (x\otimes v\otimes y,d_2(z\otimes t))_\phi&=(x\otimes v\otimes y,\sum\limits_{i=1}^n\epsilon_{a_i}za_i\otimes a_i^*\otimes t+\sum\limits_{i=1}^n\epsilon_{a_i}z\otimes a_i\otimes a_i^*t\\
 &\qquad -zb\otimes b\otimes t-z\otimes b\otimes bt)_\phi,
\end{align*}
which gives us $\forall a\in\bar Q'$
\[
(x\otimes a\otimes y,d_2(z\otimes t))_\phi=-\Tr(x\phi(t))\Tr(yza)+\Tr(x\phi(at))\Tr(yz)
\]
and
\[
 (x\otimes b\otimes y,d_2(z\otimes t))_\phi=-\Tr(x\phi(t))\Tr(yzb)-\Tr(x\phi(bt))\Tr(yz),
\]
i.e. for $v\in V$, 
\begin{align*}
 (x\otimes v\otimes y,d_2(z\otimes t))_\phi&=-\Tr(x\phi(t))\Tr(yzv)+\Tr(x\phi(vt))\Tr(yz)\\&=\Tr(xv\phi(t))\Tr(yz)-\Tr(x\phi(t))\Tr(vyz)\\
 &=(xv\otimes y-x\otimes vy,z\otimes t)_\phi=(d_1(x\otimes v\otimes y),z\otimes y)_\phi,
\end{align*}
so $d_2=d_1^\star$.
\end{proof}

\begin{proposition}
 $S_\bullet$ is an exact sequence.
\end{proposition}

\begin{proof}
 
We recall the definition of Anick's resolution \cite{An}. Denote $T_RW$ to be the tensor algebra of a graded $R$-bimodule $W$, $T_R^+W$ its augmentation ideal. Let $L\subset T_R^+W$ be an $R$-graded bimodule and $A'=T_RW/(L)$. Then we have the following resolution:

\begin{equation}
 A'\otimes_RL\otimes_RA'\stackrel{\partial}{\rightarrow}A'\otimes_RW\otimes_RA'\stackrel{f}{\rightarrow}A'\otimes_RA'\stackrel{m}{\rightarrow}A'\rightarrow 0,
\end{equation}

where $m$ is the multiplication map, $f$ is given by
\[
 f(a_1'\otimes w\otimes a_2')=a_1'w\otimes a_2'-a_1'\otimes wa_2'
\]
and $\partial$ is given by
\[\partial(a_1'\otimes l\otimes a_2')=a_1'\cdot D(l)\cdot a_2',\]
\begin{eqnarray*}
 D:T_R^+W&\rightarrow&A'\otimes_RW\otimes_RA',\\
 w_1\otimes\ldots\otimes w_n&\mapsto&\sum\limits_{p=1}^n(\overline{w_1\otimes\ldots\otimes w_{p-1}})\otimes w_p\otimes(\overline{w_{p+1}\otimes\ldots\otimes w_{n}}),
\end{eqnarray*}
where bar stands for the image in $B$ of the projection map.

In our setting, $W=V$, $L$ the $R$-bimodule generated by $\sum\limits_{i=1}^n\epsilon_{a_i}a_ia_i^*-b^2$. Then $A'=A$.\\
It is also clear that Im$(\partial)$=Im$(d_2)\subset A\otimes_RV\otimes_RA$, so from Anick's resolution we know that the part
\[A\otimes_R A[2]\stackrel{d_2}{\rightarrow}A\otimes_RV\otimes_RA\stackrel{d_1}{\rightarrow}A\otimes_RA\stackrel{d_0}{\rightarrow}A\rightarrow 0\]
is exact. Exactness of the whole complex $S_\bullet$ follows from its self duality.
\end{proof}

Since $\phi^2=1$, we can make a canonical identification
$A={}_1A_\phi\otimes_A {}_1A_\phi$ (via $x\mapsto x\otimes
1$), so by tensoring $S_\bullet$  with
${}_1A_\phi$, we obtain the exact sequence

\[
0\rightarrow A[2h]\stackrel{i'}{\rightarrow}A\otimes_R {}_1A_\phi[h+2]\stackrel{d_5}{\rightarrow}A\otimes_RV\otimes_R{}_1A_\phi[h]\stackrel{d_4}{\rightarrow}A\otimes_R{}_1A_\phi[h]\stackrel{j}{\rightarrow}{}_1A_\phi[h]\rightarrow 0.
\]
By connecting this sequence to $S_\bullet$ with $d_3=ij$ and repeating this process, we obtain the periodic Schofield resolution with period $6$:
\begin{align*}
\ldots&\rightarrow A\otimes A[2h]\stackrel{d_6}{\rightarrow}A\otimes_R {}_1A_\phi[h+2]\stackrel{d_5}{\rightarrow}A\otimes_RV\otimes_R{}_1A_\phi[h]\stackrel{d_4}{\rightarrow}A\otimes_R\mathfrak{N}[h]\\
&\stackrel{d_3}{\rightarrow}A\otimes_RA[2]\stackrel{d_2}{\rightarrow}A\otimes_RV\otimes_RA\stackrel{d_1}{\rightarrow}A\otimes_RA\stackrel{d_0}{\rightarrow}A\rightarrow0.
\end{align*}
This implies that the Hochschild homology and cohomology of $A$
is periodic with period $6$, in the sense that the shift of the
(co)homological degree by $6$ results in the shift of degree by $2h$ (respectively
$-2h$).

From that we get the periodicities for the Hochschild homology/cohomology

\begin{equation}\label{periodicity}
 HH_{j+6i}(A)\cong HH_j(A)[2ih],\qquad HH^{j+6i}(A)\cong HH^j(A)[-2ih],\qquad\forall j\geq1.
\end{equation}

\subsection{Calabi-Yau Frobenius algebras}
Let us define the functor
\[
\begin{array}{rcl}
 \mathrm{Hom}_{A^e}(-,A\otimes_\mathbb{C} A):A^e-\mathrm{mod}&\rightarrow&A^e-\mathrm{mod},\\
 M&\mapsto&M^\vee.
 \end{array}
\]
We recall the definition of the Calabi-Yau algebras from \cite{SE}.
\begin{definition}
A Frobenius algebra $A$ is called \textbf{Calabi-Yau Frobenius of dimension} \emph{m} if 
\begin{equation}
 A^\vee\simeq\Omega^{m+1}A
\end{equation}
If there is more than one such $m$, then we pick the smallest one.

If additionally $A$ has a grading, such that the above isomorphism is a graded isomorphism when composed with some shift, then we say that $A$ is a \textbf{graded Calabi-Yau Frobenius algebra}. More precisely, if $A^\vee[m']\simeq\Omega^{m+1}A$ is a graded isomorphism, where $[l]$ is the shift by $l$ with the new grading, then one says that $A$ is \textbf{graded Calabi-Yau Frobenius with dimension} $m$ \textbf{of shift} $m'$.
\end{definition}

\begin{proposition}
 $A$ is Calabi-Yau Frobenius with dimension $5$ of shift $h+2$, i.e.
 \begin{equation}
  A^\vee[h+2]\simeq\Omega^6A.
 \end{equation}
\end{proposition}
\begin{proof}
This follows from \cite{SE} since $A$ is symmetric and periodic with period $6$.
\end{proof}

From \cite{SE}, we can deduce the dualities

\begin{eqnarray}
 HH_i(A)&\cong&HH_{5-i}(A)^*[2h],\label{dual ho-ho*}\\
 HH^i(A)&\cong&HH_{5-i}(A)[-h-2],\label{dual ho-coho}\\
 HH^i(A)&\cong&HH^{11-i}(A)^*[-2h-4]\cong HH^{5-i}(A)^*[-4].
\end{eqnarray}

\subsection{Hochschild homology of $A$}\label{Hochschild homology}
Let $A^{op}$ be the algebra $A$ with opposite multiplication. We define $A^e=A\otimes_R A^{op}$. Then any $A-$bimodule naturally becomes a left $A^e-$ module (and vice versa).

Now, we apply to the Schofield resolution the functor $-\otimes_{A_e}A$ to get the  Hochschild homology complex

\begin{align*}
\ldots&\rightarrow 
A^R[2h]\stackrel{d_6'}{\rightarrow}{}_1A_\phi^R[h+2]\stackrel{d_5'}{\rightarrow}(V\otimes_R{}_1A_\phi)^R[h]\stackrel{d_4'}{\rightarrow}\\
&\stackrel{d_4'}{\rightarrow}{}_1A_\phi^R[h]\stackrel{d_3'}{\rightarrow}A^R[2]\stackrel{d_2'}{\rightarrow}(V\otimes_RA)^R\stackrel{d_1'}{\rightarrow}A^R\rightarrow0.
\end{align*}

Let $\overline{HH_i(A)}$ be $HH_0/R$ for $i=0$ and $HH_i(A)$ otherwise. We have the Connes exact sequence
\begin{equation}\label{Connes}
0\rightarrow \overline{HH_0}(A)\stackrel{B_0}{\rightarrow}\overline{HH_1}(A)\stackrel{B_1}{\rightarrow}\overline{HH_2}(A)\stackrel{B_2}{\rightarrow}\overline{HH_3}(A)\stackrel{B_3}{\rightarrow}\overline{HH_4}(A)\rightarrow\ldots,
\end{equation}
where the $B_i$ are the Connes differentials (see \cite[2.1.7.]{Lo}) and the $B_i$ are all degree-preserving.

In our case, $\overline{HH_0}(A)=A/([A,A]+R)=U^*[h]$ (see Proposition \ref{A/[A,A]}), where $U^*[h]=\langle b^i|i\,\mathrm{odd}\rangle$. From (\ref{Connes}) we know that $U^*[h]\subset\overline{HH_1(A)}$. Denote $X=\overline{HH_1(A)}/U^*[h]$. Since $\deg HH_2(A)\leq 2h$, $HH_3(A)=HH_2(A)^*[2h]$ and the Connes differential maps $HH_2(A)/X$ isomorphically to its image in $HH_3(A)$, $HH_2(A)/X$ sits in degree $h$. We call this space $K^*[h]$, where $K^*$ sits in degree $0$. $HH_3(A)=X^*[2h]\oplus K[h]$ and $HH_4(A)=U[2h]\oplus X^*[2h]$ follow from the duality (\ref{dual ho-ho*}). 
The Connes differential maps $HH_5(A)/U[2h]$ isomrphically into its image in $HH_6(A)$. Since $\deg HH_5(A)\leq 2h$ and $HH_6(A)=HH_5(A)^*[4h]$ (\ref{dual ho-ho*}), $HH_5(A)/U[2h]$ sits in degree $2h$. We call that space $Y[2h]$ where $Y$ sits in degree $0$.

From our discussion, we get the homology spaces

\begin{eqnarray*}
 HH_0(A)&=&U^*[h]\oplus R,\\
 HH_1(A)&=&U^*[h]\oplus X,\\
 HH_2(A)&=&K^*[h]\oplus X,\\
 HH_3(A)&=&K[h]\oplus X^*[2h],\\
 HH_4(A)&=&U[h]\oplus X^*[2h],\\
 HH_5(A)&=&U[h]\oplus Y[3h],\\
 HH_6(A)&=&U^*[3h]\oplus Y^*[3h],\\
 HH_{6k+i}(A)&=&HH_i(A)[2kh]\quad\forall i\geq1.
\end{eqnarray*}

\subsection{Hochschild cohomology of $A$}

We make the identifications\\
$\text{Hom}_{A^e}(A\otimes_R A,A)=A^R=\text{Hom}_{A^e}(A\otimes_R {}_1A_\phi,A)$ 
by identifying $\varphi$ with the image $\varphi(1\otimes1)=a$ (we write $\varphi=a\circ -$) and\\
$\text{Hom}_{A^e}(A\otimes_R V\otimes_R A,A)=(V\otimes_R A)^R[-2]=\text{Hom}_{A^e}(A\otimes_R V\otimes_R {}_1A_\phi,A)$ 
by identifying $\varphi$ which maps $1\otimes a\otimes1\mapsto x_a$ ($a\in\bar Q$) with the element $\sum\limits_{i=1}^n\epsilon_{a_i^*}a_i^*\otimes x_{a_i}$ (we write $\varphi=(\sum\limits_{i=1}^n\epsilon_{a_i^*}a_i^*\otimes x_{a_i}+b\otimes x_b)\circ -$).

Now, apply the functor $\text{Hom}_{A^e}(-,A)$ to the Schofield resolution to obtain the Hochschild cohomology complex

\begin{align*}
\stackrel{d_4^*}{\leftarrow}A^R[-h]\stackrel{d_3^*}{\leftarrow}A^R[-2]\stackrel{d_2^*}{\leftarrow}(V\otimes A)^R[-2]\stackrel{d_1^*}{\leftarrow}A^R&\leftarrow0
\\
\ldots\leftarrow A^R[-2h]\stackrel{d_6^*}{\leftarrow}A^R[-h-2]\stackrel{d_5^*}{\leftarrow}(V\otimes A)^R[-h-2]&\stackrel{d_4^*}{\leftarrow}
\end{align*}

and compute the differentials. We have

\[
d_1^*(x)(1\otimes y\otimes1)=x\circ d_1(1\otimes y\otimes1)=x\circ(y\otimes1-1\otimes y)=[y,x],
\]
so \[d_1^*(x)=\sum\limits_{i=1}^n\epsilon_{a_i^*}a_i^*\otimes[a_i,x]+b\otimes [b,x].\]
We have
\begin{align*}
d_2^*(\sum\limits_{i=1}^na_i\otimes x_{a_i}+b\otimes x_b)(1\otimes1)&=(\sum\limits_{i=1}^na_i\otimes x_{a_i}+b\otimes x_b)\circ(\sum\limits_{j=1}^n\epsilon_{a_j}a_j\otimes a_j^*\otimes1\\&\quad+\sum\limits_{j=1}^n\epsilon_{a_j}1\otimes a_j\otimes a_j^*-b\otimes b\otimes 1-1\otimes b\otimes b)\\
&=\sum\limits_{i=1}^n(a_ix_{a_i}-x_{a_i}a_i)-(bx_b+x_bb)\\&=\sum\limits_{i=1}[a_i,x_{a_i}]-(bx_b+x_bb),
\end{align*}
so
\[
d_2^*(\sum\limits_{i=1}^na_i\otimes x_{a_i}+b\otimes x_b)=\sum\limits_{i=1}^n[a_i,x_{a_i}]-(bx_b+x_bb).
\]
We have
\[
d_3^*(x)(1\otimes1)=x\circ d_3(1\otimes1)=x\circ(\sum\limits_{x_i\in\mathfrak{B}} \phi(x_i)\otimes x_i^*)=\sum\limits_{x_i\in\mathfrak{B}} \phi(x_i)x x_i^*=0,
\]

so
\[
d_3^*(x)=\sum\limits_{x_i\in\mathfrak{B}} \phi(x_i)xx_i^*,
\]
and we evaluate this sum:
let $\rho=\sum\limits_{i=1}^{n}(-1)^{n+i}e_i$, then $\rho^2=\sum\limits_{i=1}^{n}e_i=1$, and for all monomials $x\in A$, $\rho\phi(x)\rho=(-1)^{\deg(x)}x$. We write $x=\rho y$ (where $y=\rho x$), then

\[
 d_3^*(x)=\sum\limits_{x_i\in\mathfrak{B}} \phi(x_i)xx_i^*=\rho\sum\limits_{x_i\in\mathfrak{B}} (\rho\phi(x_i)\rho)yx_i^*=\rho\sum\limits_{x_i\in\mathfrak{B}} (-1)^{\deg(x_i)}x_iyx_i^*.
\]

The map $y\mapsto \sum\limits_{x_i\in\mathfrak{B}} (-1)^{\deg(x_i)}x_iyx_i^*$ is zero in positive degree, and the restriction to $\deg y=0$ is a map $\bigoplus_{i\in I}\mathbb{C}e_i\rightarrow\bigoplus_{i\in I}\mathbb{C}\omega_i$, given by the matrix $H_A(-1)=(1+(-1)^{h})(2+C)=0$, since the Coxeter number $h=2n+1$ is odd. This implies that 

\[
 d_3^*=0.
\]
We have
\[
d_4^*(x)(1\otimes a_i\otimes1)=x\circ d_1(1\otimes a_i\otimes1)=x\circ(a_i\otimes1-1\otimes a_i)=a_ix-x\phi(a_i)=[a,x],
\]
and
\[
d_4^*(x)(1\otimes b\otimes 1)=x\circ d_1(1\otimes b\otimes1)=x\circ(b\otimes1-1\otimes b)=bx-x\phi(b),
\]

so \[d_4^*(x)=\sum\limits_{i=1}^n\epsilon_{a_i^*}a_i^*\otimes [a_i,x]+b\otimes(xb+bx),\]
We have
\begin{align*}
d_5^*(\sum\limits_{i=1}^na_i\otimes x_{a_i}+b\otimes x_b)(1\otimes1)&=(\sum\limits_{i=1}^na_i\otimes x_{a_i}+b\otimes x_b)\circ(\sum\limits_{i=1}^n\epsilon_{a_i}a_i\otimes a_i^*\otimes1\\&\quad+\sum\limits_{i=1}^n\epsilon_{a_i}1\otimes a_i\otimes a_i^*)\\
&=\sum\limits_{i=1}^n(a_ix_{a_i}-x_{a_i}\phi(a_i))+(bx_b+x_b\phi(b)),
\end{align*}
so
\[
d_5^*(x,y)=[x,y].
\]
We have
\[
d_6^*(x)(1\otimes1)=x\circ d_6(1\otimes1)=x\circ(\sum\limits_{x_i\in\mathfrak{B}} \phi(x_i)\otimes x_i^*)=\sum\limits_{x_i\in\mathfrak{B}} \phi(x_i)x\phi(x_i^*)=-\sum\limits_{x_i\in\mathfrak{B}} x_ixx_i^*,
\]
so
\[
d_6^*(x)=-\sum\limits_{x_i\in\mathfrak{B}} x_ixx_i^*.
\]

From our results about Hochschild homology and the dualities (\ref{dual ho-coho}), we obtain the following spaces for the Hochschild cohomology (for $HH^0(A)$, keep in mind that we get $HH^6(A)=HH^0(A)[-2h]/\text{Im}d_6^*$, and the image of $d_6^*$ lies in top degree).
The cohomology spaces are
\begin{eqnarray*}
 HH^0(A)&=&U[-2]\oplus L[h-2],\\
 HH^1(A)&=&U[-2]\oplus X^*[h-2],\\
 HH^2(A)&=&K[-2]\oplus X^*[h-2],\\
 HH^3(A)&=&K^*[-2]\oplus X[-h-2],\\
 HH^4(A)&=&U^*[-2]\oplus X[-h-2],\\
 HH^5(A)&=&U^*[-2]\oplus Y^*[-h-2],\\
 HH^6(A)&=&U[-2h-2]\oplus Y[-h-2],\\
 HH^{6k+i}(A)&=&HH^i(A)[-2kh]\quad\forall{i\geq 1}.
\end{eqnarray*}

We have $L[h-2]=R^*[h-2]$. Since there is no non-top degree element in $A$ which commutes with all $a\in\bar Q'$ and anticommutes with $b$, $\ker d_4^*$ lies in top degree $-2$ which implies that the space $X$ has to be zero.

From the discussion in Subsection \ref{Hochschild homology}, we know that $K$ is a  degree-zero space, so $HH^2(A)$ sits entirely in degree $-2$. Since $d_3^*=0$ and the image of $d_4^*$ lies in degree $>-2$, $K=\bigoplus_{i\in I}\mathbb{C} e_i$, so $K$ is $n$-dimensional. This proves Theorem \ref{UK} (b).

The map $d_6^*$ can be viewed as a map $\bigoplus_{i\in I}\mathbb{C}e_i\rightarrow\bigoplus_{i\in I}\mathbb{C}\omega_i$, given by the matrix $-H_A(1)=-2(2-C)^{-1}$. Since it is nondegenerate, the space $Y$ is also zero.

Theorems \ref{cohomology} and \ref{homology} follow.

\subsection{Cyclic homology of $A$}
Given the Connes exact sequence (\ref{Connes}), we define
\begin{align*}
\overline{HC_i}(A)&=\ker(B_{i+1}:\overline{HH_{i+1}}(A)\rightarrow\overline{HH_{i+2}}(A))\\
&=\mbox{Im}(B_i:\overline{HH_i}(A)\rightarrow\overline{HH_{i+1}}(A)).
\end{align*}

The usual cyclic homology $HC_i(A)$ is related to the reduced one
by the equality $\overline{HC}_i(A)=HC_i(A)$ for $i>0$, 
and $\overline{HC}_0(A)=HC_0(A)/R$.

We write down Connes exact sequence, together with the Hochschild homology spaces:
$$
\CD
                      0                    \\
                    @VVV\\
       \overline{HH_0}(A)@=U^*[h]     \\
                   @VB_1VV    @V\sim VV\\		    
       \overline{HH_1}(A)@=   U^*[h]@.  \overline{HC}_0(A)=U^*[h]\\
                   @VB_1VV    @V 0 VV\\
       \overline{HH_2}(A)@=   K^*[h]@.  HC_1(A)=0\\
                   @VB_2VV            @V\sim VV\\
       \overline{HH_3}(A)@= K[h]   @.   HC_2(A)=K^*[h] \\
                   @VB_3VV     @V 0 VV\\
       \overline{HH_4}(A)@= U[h] @.     HC_3(A)=0\\
                   @VB_4VV       @V\sim VV \\
       \overline{HH_5}(A)@=   U[h]     @.     HC_4(A)=U[h]\\
                   @VB_5VV     @V 0 VV\\
       \overline{HH_6}(A)@=   U^*[3h]     @.  HC_5(A)=0\\
                   @VB_6VV       @V0VV    \\
       \overline{HH_7}(A)@=U^*[3h]        @.  HC_6(A)=U^*[3h]\\
                           @VB_7VV\\
                              \vdots
\endCD
$$
This proves Theorem \ref{cyclic}.

\section{Basis of $HH^*(A)$}\label{basis}
Now we construct a basis of $HH^*(A)$.
\subsection{$HH^0(A)=Z$}
We compute the structure of the center.

\begin{proposition}
 The non-topdegree central elements lie in even degrees, one in each degree (up to scaling). They are given by
  \begin{equation}\
  z_{2k}=\sum\limits_{i=k+1}^n c_{i,2k},\qquad 0\leq k\leq\frac{h-3}{2}
 \end{equation}
\end{proposition}
\begin{proof}
 First we prove that a degree $2k$-element $z$ is a multiple of $z_{2k}$:
 $z$ commutes with all $e_i$, hence lies in $\bigoplus_{i\in A}e_iAe_i$. From the discussion, in Subsection (\ref{bases}), we can write \[z=\sum\limits_{i=k+1}^n\lambda_ic_{i,2k}=\sum\limits_{i=k+1}^n\lambda_i(a_i^*a_i)^k.\] 
 Now,
 \[
  \forall j\geq k+1,\quad\lambda_j\underbrace{a_j(a_j^*a_j)^k}_{\neq0}=a_jz=za_j=\lambda_{j+1}(a_ja_j^*)^ka_j
 \]
which implies that all $\lambda_i$ are equal. So each even degree central is a multiple of $z_{2k}$. 

Since $z_{2k}^*=z_{2k}$ and $z_{2k}$ commutes with all $a_j$, $z_{2k}$ also commutes with all $a_j^*$. Commutativity with $b$ is clear, since each element in $e_nAe_n$ can be expressed as a polynomial in $b$.

So $z_{2k}=\sum\limits_{i=k+1}^n c_{i,2k}$ is the central element in degree $2k$.

Now, let $z$ be of odd degree $<h-2$. From Subsection (\ref{bases}), we can write
\[
 z=\sum\limits_{i=n-k}^n\lambda_ic_{i,2k+1},
\]
where 
\[
c_{i,2k+1}=a_i^*\cdot\ldots\cdot a_{n-1}^*b^{2i-2n+2k}a_{n-1}\ldots\cdot a_i
\]
and the law $\forall i\geq n-k$,
\begin{align*}
 (0\neq)a_ic_{i,2k+1}&=a_{i+1}^*\cdot\ldots\cdot a_{n-1}^*b^{2i-2n+2k+2}a_{n-1}\ldots\cdot a_i=c_{i+1,2k+1}a_i,
\end{align*}

so we have 
\[
\lambda_{n-k}\underbrace{c_{n-k,2k+1}a_{n-k-1}}_{\neq0}=za_{n-k-1}=a_{n-k-1}z=0
\]
and 
\[
 \forall j\geq n-k,\quad \lambda_ja_jc_{j,2k+1}=a_jz=za_j=\lambda_{j+1}c_{i+1,2k+1}a_j=\lambda a_jc_{j,2k+1}.
\]
This implies that all $\lambda_j=0$, so we have no non-top odd degree central elements.

\end{proof}

Theorem \ref{UK} (a) follows.

\subsection{$HH^1(A)$}
Since $HH^1(A)=U[-2]$, we know from the previous subsection that the Hilbert series of $HH^1(A)$ is $\sum\limits_{i=0}^{\frac{h-3}{2}}t^{2i}$. It is easy to see that 
\[
\theta_{2k}:=\sum\limits_{i=1}^na_i\otimes a_i^*z_{2k}-\sum\limits_{i=1}^na_i^*\otimes a_iz_{2k}+b\otimes bz_{2k},\qquad 0\leq k\leq \frac{h-3}{2}
\]
lie in $\ker d_1^*$. The cup product calculation $HH^1(A)\cup HH^4(A)$ will show that each $\theta_{2k}$ is nonzero (since the product with $\zeta_{2k}$ is nonzero).

\subsection{$HH^2(A)$ and $HH^3(A)$}
$HH^2(A)$ and $HH^3(A)$ sit in degree $-2$ and both are $n$-dimensional. So $HH^2(A)$ is the bottom degree part of $A^R[-2]$ and $HH^3(A)$ the top degree part of $A^R[-h]$. Denote $f_i=[e_i]\in HH^2(A)$ and $h_i=[\omega_i]\in HH^3(A)$, we have 
\[
 HH^2(A)=\bigoplus\limits_{i=1}^n\mathbb{C}f_i,\qquad HH^3(A)=\bigoplus\limits_{i=1}^n\mathbb{C}h_i.
\]

\subsection{$HH^4(A)$}
The Hilbert series of $HH^4(A)$ is $t^{-4}\sum\limits_{i=0}^{\frac{h-3}{2}}t^{-2i}$. We claim that a basis is given by
\[
 \zeta_{2i}:=[-b\otimes b^{h-3-2i}].
\]
It is clear that $\zeta_{2i}$ all lie in $\ker d_5^*$. Since the image of $d_4^*$ has zero trace, $\zeta_0$ is nonzero in $HH^4(A)$. And $\zeta_{2i}\neq0$ follows from $z_{2i}\zeta_{2i}=\zeta_0$.

\subsection{$HH^5(A)$} From Proposition \ref{A/[A,A]}, we know that the space\\
$HH^5(A)=A^R/([A,A]^R+R)[-h-2]=A/([A,A]+R)[-h-2]$ is spanned by 
\[
 \psi_{2i}:=[b^{h-3-2i}].
\]

\section{The Hochschild cohomology ring $HH^*(A)$}
The degree ranges of the Hochschild cohomology spaces are
\begin{eqnarray*}
 0&\leq\deg HH^0(A)&\leq h-2,\\
 0&\leq\deg HH^1(A)&\leq h-3,\\
-2&=\deg HH^2(A),\\
-2&=\deg HH^3(A),\\
 -h-1&\leq\deg HH^4(A)&\leq -4,\\
 -h-1&\leq\deg HH^5(A)&\leq -4,\\
 -2h&\leq\deg HH^6(A)&\leq -h-3.
\end{eqnarray*}

We compute the cup product in terms of our constructed basis in $HH^*(A)$ from the last section.

\subsection{The $Z$-module structure of $HH^*(A)$}
$HH^0(A)$ is a local ring, with radical generated by $z_2$.
In $HH^0(A)$, we have $z_{2i}z_{2j}=z_{2(i+j)}$ for $2i+2j\leq h-3$, and the product is $0$ otherwise. $HH^i(A)$ are cyclic $Z$-modules for $i=1,4,5$, generated by $\theta_0,\,\zeta_{h-3},\,\psi_{h-3}$ respectively. The $Z$-modules $HH^2(A)$ and $HH^3(A)$ are annihilated by the radical of $Z$.

\subsection{$HH^i(A)\cup HH^j(A)$ for $i,j$ odd}
All cup products $HH^i(A)$ with $HH^j(A)$ for $i,j$ odd are zero by degree argument.

\subsection{$HH^1(A)\cup HH^2(A)$}
By degree argument, $\theta_i f_j=0$ for $i\neq0$. 
\begin{proposition}\label{H1xH2}
 The multiplication with $\theta_0$ gives us a map
 \[HH^2(A)=K[-2]\stackrel{\alpha}{\rightarrow} K^*[-2]=HH^3(A),\]
 given by the matrix
 \[
 h
 \left[
 \begin{array}{cccccc}
 2&-1&0&\cdots&\cdots&0\\
 -1&2&\ddots&\ddots&&\vdots\\
 0&\ddots&\ddots&\ddots&\ddots&\vdots\\
 \vdots&\ddots&\ddots&\ddots&\ddots&0\\
 0&\cdots&0&-1&2&-1\\
 0&\cdots&\cdots&0&-1&3
 \end{array}
 \right]^{-1}.\]
\end{proposition}

\begin{proof}

Let $x\in K[-2]$, represented by the map
 \begin{eqnarray*}
  f_x:A\otimes A[2]&\longrightarrow& A,\\
  1\otimes1&\longmapsto&x,
 \end{eqnarray*}
which we lift to
 \begin{eqnarray*}
  \hat f_x:A\otimes A[2]&\longrightarrow& A\otimes A,\\
  1\otimes1&\longmapsto&1\otimes x.
 \end{eqnarray*}
 Then we have
\[
 \hat f_xd_3(1\otimes 1)=\hat f_x(\sum\limits_{x_j\in\mathfrak{B}} \phi(x_j)\otimes x_j^*)=\sum\limits_{x_j\in \mathfrak{B}} \phi(x_j)\otimes xx_j^*.
\]
To compute the lift $\Omega f_i$, we need to find out the preimage of $\sum\limits_{x_j\in\mathfrak{B}} \phi(x_j)\otimes xx_j^*$ under $d_1$.

\begin{definition}
 Let $b_1,\ldots, b_k$ be arrows, $p$ the monomial $b_1\cdots b_k$ and define 
 \[
v_p:=(1\otimes b_1\otimes b_2\cdots b_k+b_1\otimes b_2\otimes b_3\cdots b_k+\ldots+b_1\cdots b_{k-1}\otimes b_k\otimes 1).
 \]
\end{definition}
We will use the following lemma in our computations.
\begin{lemma}
In the above setting,
 \[
  d_1(v_p)=(b_1\cdots b_k\otimes 1-1\otimes b_1\cdots b_k).
 \]
\end{lemma}
From that, we immediately see that
\[
 \sum\limits_{x_j\in\mathfrak{B}}\phi(x_j)\otimes xx_j^*=d_1(\sum\limits_{x_j\in \mathfrak{B}}v_{\phi(x_j)}xx_j^*)+1\otimes \underbrace{\sum\limits_{x_j\in \mathfrak{B}}\phi(x_j)xx_j^*}_{=0},
\]
so we have
\begin{eqnarray*}
 \Omega f_x:\Omega^3(A)&\rightarrow&\Omega(A),\\
 1&\mapsto&\sum\limits_{x_j\in \mathfrak{B}}v_{\phi(x_j)}xx_j^*.
\end{eqnarray*}
Then we have
\[
 \theta_0(\sum\limits_{x_j\in \mathfrak{B}}v_{\phi(x_j)}xx_j^*)=\sum\limits_{x_j\in \mathfrak{B}}\deg(x_j)\phi(x_j)xx_j^*.
\]
So we get
\[
 \theta_0f_x=\sum\limits_{x_j\in \mathfrak{B}}\deg(x_j)\phi(x_j)xx_j^*=\sum\limits_{k,l=1}^n\sum\limits_{x_j\in \mathfrak{B}_{kl}}\deg(x_j)\phi(x_j)xx_j^*.
\]
$\phi(x_j)$ is $x_j$ if the number of $b's$ in $x_j$ is even and $-x_j$ if it is odd. Observe that the number of $b's$ in $x_j$ and $\deg(x_j)-d(k,l)$ (where $d(k,l)$ is the distance between the vertices $k$ and $l$) have the same parity. So $\phi(x_j)=(-1)^{\deg(x_j)-d(k,l)}x_j$, and so the multiplication with $\theta_0$ induces a map
\begin{equation}\label{alpha}
 HH^2(A)=K[-2]\stackrel{\alpha}{\rightarrow} K^*[-2]=HH^3(A),
\end{equation}
given by the matrix 
\begin{equation}
 \left(H_A^\phi\right)_{k,l}=\sum\limits_{h_j\in \mathfrak{B}_{k,l}}(-1)^{\deg(x_j)-d(k,l)}\deg(x_j)=(-1)^{d(i,j)}\left.\left(\frac{d}{dt}H_A(t)_{k,l}\right)\right|_{t=-1}.
\end{equation}
Let us define
\[H_A^\delta:=\left.\left(\frac{d}{dt}H_A(t)\right)\right|_{t=-1}.\]

Then we have
\begin{align*}
 H_A^\delta&=\left.\left((1+t^h)\frac{d}{dt}(1-Ct+t^2)^{-1}+ht^{h-1}(1-Ct+t^2)^{-1}\right)\right|_{t=-1}\\
 &=h(2+C)^{-1}.
\end{align*}

For any nondegenerate matrix $M$, call $M_-$ the matrix obtained from $M$ by changing all signs in the $(i,j)$-entry whenever $d(i,j)$ is odd. It is easy to see that for matrices $M=N^{-1}$, $M_-=(N_-)^{-1}$. In our case, we have $H_A^\phi=(H_A^\delta)_-$. This implies 
\begin{equation}
 H_A^\phi=h((2+C)_-)^{-1}.
\end{equation}
\end{proof}
\subsection{$HH^1(A)\cup HH^4(A)$}\label{H1xH4}
Since $HH^1(A)=Z\theta_0$ and $HH^4(A)=Z\zeta_{h-3}$, it is enough to compute $\theta_0\zeta_{h-3}$.

\begin{proposition}
 Given $\theta_0\in HH^1(A)$ and $\zeta_{h-3}\in HH^4(A)$, we get the cup product
 \begin{equation}\label{thetaxzeta}
  \theta_0\zeta_{h-3}=\psi_{h-3}.
 \end{equation}
\end{proposition}
\begin{proof}
$\zeta_{h-3}$ represents the map

\begin{eqnarray*}
 \zeta_{h-3}:A\otimes V\otimes {}_1A_\phi[h]&\rightarrow&A,\\
 1\otimes b\otimes 1&\mapsto& -e_n,\\
 1\otimes a_i\otimes 1&\mapsto&0,\\
 1\otimes a_i^*\otimes 1&\mapsto&0,
\end{eqnarray*}
and it lifts to

\begin{eqnarray*}
 \hat\zeta_{h-3}:A\otimes V\otimes {}_1A_\phi[h]&\rightarrow&A\otimes A,\\
 1\otimes b\otimes 1&\mapsto& -e_n\otimes e_n,\\
 1\otimes a_i\otimes 1&\mapsto&0,\\
 1\otimes a_i^*\otimes 1&\mapsto&0.
\end{eqnarray*}
Then
\begin{eqnarray*}
 (\hat\zeta_{h-3}\circ d_5)(1\otimes 1)&=&\hat\zeta_{h-3}(\sum\limits_{i=1}\epsilon_{a_i}a_i\otimes a_i^*\otimes1+\sum\limits_{i=1}\epsilon_{a_i}\otimes a_i\otimes a_i^*\\&&-b\otimes b\otimes 1-1\otimes b\otimes b) \\&=& b\otimes1-1\otimes b=d_1(1\otimes b\otimes 1),
\end{eqnarray*}

so we have 

\begin{eqnarray*}
 \Omega\zeta_{h-3}:\Omega^5(A)&\rightarrow&\Omega(A),\\
 1\otimes 1&\mapsto&1\otimes b\otimes 1,
\end{eqnarray*}
and this gives us
\[
 (\theta_0\circ\zeta_{h-3})(1\otimes1)=b,
\]
so the cup product is
\begin{equation}
 \theta_0\zeta_{h-3}=[b]=\psi_{h-3}.
\end{equation}
\end{proof}

\subsection{$HH^2(A)\cup HH^3(A)$}\label{H2xH3}
We compute the cup product in the following proposition.
\begin{proposition}\label{perfect}
 For the basis elements $f_i\in HH^2(A)$, $h_j\in HH^3(A)$, the cup product is
 \begin{equation}
  f_ih_j=\delta_{ij}\psi_0.
 \end{equation}
\end{proposition}
\begin{proof}
 Recall the maps 
 \begin{eqnarray*}
 h_j:A\otimes {}_1A_{\phi}&\rightarrow&A,\\
 1\otimes1&\mapsto&\omega_j
 \end{eqnarray*}
 
 and lift them to 
 \begin{eqnarray*}
 \hat h_j:A\otimes {}_1A_{\phi}&\rightarrow&A\otimes A,\\
 1\otimes1&\mapsto&1\otimes\omega_j.
 \end{eqnarray*}
 Then $\forall\, a\in\bar Q$ we have 
 \[
 \hat h_j(d_4(1\otimes a\otimes 1))=\hat h_j(a\otimes1-1\otimes a)=a\otimes\omega_j=d_1(1\otimes a\otimes\omega_j),
 \]
so
\begin{eqnarray*}
 \Omega h_j:\Omega^4(A)&\rightarrow&\Omega(A),\\
 1\otimes a\otimes 1&\mapsto&1\otimes a\otimes\omega_j.
\end{eqnarray*}
Then we have
\begin{eqnarray*}
 \Omega h_j(d_5(1\otimes1))&=&\Omega h_j(\sum\limits_{i=1}^n\epsilon_{a_i}a_i\otimes a_i^*\otimes1+\sum\limits_{i=1}^n\epsilon_{a_i}1\otimes a_i\otimes a_i^*\\
 &&-b\otimes b\otimes1-1\otimes b\otimes b)\\
 &=&(\sum\limits_{i=1}^n\epsilon_{a_i}a_i\otimes a_i^*-b\otimes b\otimes\omega_j=d_2(1\otimes\omega_j),
\end{eqnarray*}
so
\begin{eqnarray*}
 \Omega^2 h_j:\Omega^5(A)&\rightarrow&\Omega^2(A),\\
 1\otimes 1&\mapsto&1\otimes \omega_j.
\end{eqnarray*}
This gives us 
\[
 f_i(\Omega^2 h_j)(1\otimes1)=f_i(1\otimes\omega_j)=\delta_{ij}\omega_j,
\]
i.e. the cup product
\[
 f_ih_j=\delta_{ij}[\omega_j]=\delta_{ij}\psi_0.
\]
\end{proof}

\subsection{$HH^2(A)\cup HH^2(A)$}\label{H2xH2}
Since $\deg HH^2(A)=-2$, their product has degree $-4$ (i.e. lies in $\text{span}(\zeta_0)$), so it can be written as

\begin{eqnarray*}
 HH^2(A)\times HH^2(A)&\rightarrow&HH^4(A),\\
 (x,y)&\mapsto&\langle-,-\rangle\zeta_0,
\end{eqnarray*}
where $\langle-,-\rangle:HH^2(A)\times HH^2(A)\rightarrow\mathbb{C}$ is a bilinear form.

\begin{proposition}
 The cup product $HH^2(A)\times HH^2(A)\rightarrow HH^4(A)$ is given by $\langle-,-\rangle=\alpha$, where $\alpha$ (from Proposition \ref{alpha}) is regarded as a symmetric bilinear form.
\end{proposition}
\begin{proof}

We use (\ref{thetaxzeta}) to get
\begin{equation}\label{H2xH2 1}
 \theta_0(f_if_j)=\theta_0(\langle f_i,f_j\rangle\zeta_0)=\langle f_i,f_j\rangle\psi_0.
\end{equation}

On the other hand, by Proposition \ref{alpha} and Proposition \ref{perfect},
\begin{equation}\label{H2xH2 2}
 (\theta_0f_i)f_j=\alpha(f_i)f_j=\sum\limits_{l=1}^n\left(H_A^\phi\right)_{li}h_lf_j=\left(H_A^\phi\right)_{ji}\psi_0=\left(H_A^\phi\right)_{ij}\psi_0.
\end{equation}
By associativity of the cup product, we can equate (\ref{H2xH2 1}) and (\ref{H2xH2 2}) to get
\begin{equation}
 \langle f_i,f_j\rangle=\left(H_A^\phi\right)_{ij}.
\end{equation}
\end{proof}

\subsection{$HH^2(A)\cup HH^4(A)$}
By degree argument, $f_i\zeta_{j}=0$ for $j< h-3$ and $f_i\zeta_{h-3}=\lambda_i \varphi_0(z_{h-3})$ for some $\lambda_i\in\mathbb{C}$. 
\begin{proposition}
 We have
 \begin{equation}\label{H2xH4}
  f_i\zeta_{h-3}=i\cdot z_{h-3}.
 \end{equation}
\end{proposition}

\begin{proof}
Let $x\in HH^2(A)$. $x$ is represented by a map $f_x$,

\begin{eqnarray*}
 f_x:A\otimes A[2]&\rightarrow&A,\\
 1\otimes1&\mapsto&x,
\end{eqnarray*}
and we lift it to 
\begin{eqnarray*}
 f_x:A\otimes A[2]&\rightarrow&A,\\
 1\otimes1&\mapsto&1\otimes x.
\end{eqnarray*}
We know that for $h_j\in HH^3(A)$ and $x=\sum\limits_{i=1}^nr_if_i$ the cup product is $xh_j=r_jf_j$. This determines the lift 
\begin{eqnarray*}
\Omega^3 f_x:\Omega^5(A)&\rightarrow&\Omega^3(A),\\
1\otimes 1&\mapsto&x.
\end{eqnarray*}

 Then 
\begin{eqnarray*}
 \Omega^4f_xd_6(1\otimes 1)&=&\Omega^4 f_x(\sum\limits_{x_j\in\mathfrak {B}} \phi(x_j)\otimes x_j^*)=\sum\limits_{x_j\in \mathfrak{B}} \phi(x_j)\otimes x\phi(x_j^*)\\
 &=&-\sum\limits_{x_j\in \mathfrak{B}} x_j\otimes xx_j^*=d_4(-\sum\limits_{x_j\in \mathfrak{B}} v_{x_j}x\phi(x_j^*)).
\end{eqnarray*}
For each term $v_{x_j}x\phi(x_j^*)$,

\begin{eqnarray*}
 \zeta_{h-3}(v_{x_j}x\phi(x_j^*))=\left\{
 \begin{array}{cl}
0&\text{if }x_j\text{ contains even number of }b's\\
-\frac{x_j}{b}xx_j^*&\text{if }x_j\text{ contains odd number of }b's,
 \end{array}
 \right.
\end{eqnarray*}
where for a monomial $x_j$, the expression ''$\frac{x_j}{b}$'' means removing one letter $b$ (and it doesn't matter which one you remove). Denote $\mathfrak{B}^{odd}$ (resp. $\mathfrak{B}^{even}$ a basis of $e_kAe_l$ which have odd (resp. even) number of $b$'s in their monomial expression. Then 
\[
 \zeta_{h-3}\circ \Omega^4f_x(1\otimes 1)=\sum\limits_{x_j\in \mathfrak{B}^{odd}}\frac{x_jxx_j^*}{b}.
\]
The automorphism $\gamma$ which reverses all arrows of a path is the identity on $A^{top}$. Let $(x_i)$ be a basis of $A$, $(x_i^*)$ its dual basis. Then $(\gamma(x_i))$ is a basis and $(\gamma(x_i^*))$ its dual basis. This shows that
\[
 \sum\limits_{x_j\in \mathfrak{B}^{odd}}x_jxx_j^*=\sum\limits_{x_j\in \mathfrak{B}^{odd}}\frac{x_j^*xx_j}{b}=\sum\limits_{x_j\in \mathfrak{B}^{even}}\frac{x_jxx_j^*}{b},
\]
so 
\[
 \sum\limits_{x_j\in \mathfrak{B}^{odd}}x_jxx_j^*=\frac{1}{2}\sum\limits_{x_j\in\mathfrak{B}}\frac{x_jxx_j^*}{b}.
\]
The $(h-3)$-degree part of $A$ lies in $e_nAe_n$ and is spanned by $z_{h-3}b^{h-3}$. This means that $\frac{\omega_n}{b}=z_{h-3}$ and $\frac{\omega_i}{b}=0$ for $i<n$. We get

\begin{equation}
 \zeta_{h-3}f_l=\frac{1}{2}\left(H_A(1)\right)_{nl}\varphi_0(z_{h-3})=l\cdot \varphi_0(z_{h-3}).
\end{equation}
\end{proof}

\subsection{$HH^2(A)\cup HH^5(A)$}
By degree argument, $f_i\psi_j=0$ for $j\neq h-3$.

\begin{proposition}
 We have
 \begin{equation}
  f_i\zeta_{h-3}=i\cdot\varphi_0(\theta_{h-3}).
 \end{equation}
\end{proposition}
\begin{proof}
 Since $\psi_{h-3}=\theta_0\zeta_{h-3}$, we have
 \[
  f_i\psi_{h-3}=(f_i\zeta_{h-3})\theta_0=i\cdot \varphi_0(z_{h-3})\theta_0=i\cdot\varphi_0(\theta_{h-3}).
 \]
\end{proof}

\subsection{$HH^3(A)\cup HH^4(A)$}
$h_i\zeta_{j}=0$ for $j< h-3$ and $h_i\zeta_{h-3}=\lambda_i \varphi_0(\theta_{h-3})$ for some $\lambda_i\in\mathbb{C}$. 
\begin{proposition}
 We have 
 \begin{equation}\label{H3xH4}
  h_i\zeta_{h-3}=\delta_{in}\varphi_0(\theta_{h-3}).
 \end{equation}
\end{proposition}

\begin{proof}
Let $\lambda_i$ be from above.
From (\ref{H2xH4}), we get
 \[
  \theta_0f_i\zeta_{h-3}=i\cdot \varphi_0(\theta_{h-3}),
 \]
and we use (\ref{H1xH2}) to see that 
\[
 \left(
 \begin{array}{c}
   \lambda_1 \\
   \lambda_2\\
   \vdots\\
   \vdots\\
   \vdots\\
   \lambda_n
 \end{array}
\right)
=\frac{1}{2n+1}
\left[
 \begin{array}{cccccc}
 2&-1&0&\cdots&\cdots&0\\
 -1&2&\ddots&\ddots&&\vdots\\
 0&\ddots&\ddots&\ddots&\ddots&\vdots\\
 \vdots&\ddots&\ddots&\ddots&\ddots&0\\
 0&\cdots&0&-1&2&-1\\
 0&\cdots&\cdots&0&-1&3
 \end{array}
 \right]
  \left(
 \begin{array}{c}
   1\\
   2\\
   \vdots\\
   \vdots\\
   \vdots\\
   n
 \end{array}
 \right)=
 \left(
 \begin{array}{c}
   0\\
   0\\
   \vdots\\
   \vdots\\
   \vdots\\
   1
 \end{array}
 \right).
\]
\end{proof}
 
\subsection{$HH^4(A)\cup HH^4(A)$}
By degree argument, $\zeta_i\zeta_j=0$ if $i<h-3$ or $j<h-3$ and $\zeta_{h-3}^2=\sum\limits_{k=1}^n\lambda_k\varphi_0(f_k)$.
\begin{proposition}
 We have
 \begin{equation}\label{H4xH4}
  \zeta_{h-3}^2=\varphi_0(f_n).
 \end{equation}
\end{proposition}
\begin{proof}
 Let $\lambda_k$ be from above. Then we have, using (\ref{perfect}),
 \[
  h_l\zeta_{h-3}^2=\lambda_l\psi_0.
 \]
Using (\ref{H3xH4}), the LHS becomes \[\delta_{ln}\theta_{h-3}\zeta_{h-3}=\delta_{ln}\psi_0,\]
so 
\[
 \lambda_{l}=\delta_{ln}.
\]
\end{proof}

\subsection{$HH^4(A)\cup HH^5(A)$}
By degree argument, $\zeta_i\psi_j=0$ if $i<h-3$ or $j<h-3$.
\begin{proposition}
 We have
 \begin{equation}
 \zeta_{h-3}\psi_{h-3}=\sum\limits_{i=1}^ni\varphi_0(h_i).
 \end{equation}
\end{proposition}
\begin{proof}
 We use (\ref{H1xH2}), (\ref{H1xH4}) and (\ref{H4xH4}) to obtain
 \[
  \zeta_{h-3}\psi_{h-3}=\zeta_{h-3}^2\theta_0=f_n\theta_0=\sum\limits_{i=1}^ni\varphi_0(h_i).
 \]
The last equality follows from
\[
  (2n+1)\left[
 \begin{array}{cccccc}
 2&-1&0&\cdots&\cdots&0\\
 -1&2&\ddots&\ddots&&\vdots\\
 0&\ddots&\ddots&\ddots&\ddots&\vdots\\
 \vdots&\ddots&\ddots&\ddots&\ddots&0\\
 0&\cdots&0&-1&2&-1\\
 0&\cdots&\cdots&0&-1&3
 \end{array}
 \right]^{-1}
  \left(
 \begin{array}{c}
   0\\
   0\\
   \vdots\\
   \vdots\\
   \vdots\\
   1
 \end{array}
 \right)=
   \left(
 \begin{array}{c}
   1\\
   2\\
   \vdots\\
   \vdots\\
   \vdots\\
   n
 \end{array}
 \right).
\]
\end{proof}

\section{Batalin-Vilkovisky structure on Hochschild cohomology}
From general theory, we have an isomorphism $\mathbb{D}: HH_\bullet(A)\rightarrow HH^{6m+5-\bullet}(A)$ $\forall m\geq0$. It translates the Connes differential $B:HH_\bullet(A)\rightarrow HH_{\bullet+1}(A)$ on Hochschild homology into a differential $\Delta:HH^\bullet(A)\rightarrow HH^{\bullet -1}(A)$ on Hochschild cohomology, i.e. we have the commutative diagram

$$
\CD
HH_\bullet(A)     @>B>>  HH_{\bullet+1}(A)\\
@V\mathbb{D}V\sim V      @V\sim V\mathbb{D}V\\         
HH^{6m+5-\bullet}(A)[(2m+1)h+2] @>\Delta>>  HH^{6m+4-\bullet}(A)[(2m+1)h+2]
\endCD
$$

\begin{theorem}\emph{(BV structure on Hochschild cohomology)}
 $\Delta$ makes $HH^\bullet(A)$ a Batalin-Vilkovisky algebra, i.e. 
 for the Gerstenhaber bracket we get the following equation:
\begin{equation}\label{BV-identity}
 [a,b]=\Delta(a\cup b)-\Delta(a)\cup b-(-1)^{|a|}a\cup\Delta(b),\qquad\forall a,b\in HH^*(A).
\end{equation}

The isomorphism $\mathbb{D}$ intertwines contraction and cup-product maps, i.e. we have
\begin{equation}\label{intertwining}
 \mathbb{D}(\iota_\eta c)=\eta\cup\mathbb{D}(c),\qquad\forall c\in HH_\bullet(A),\,\eta\in HH^\bullet(A).
\end{equation}
\end{theorem}

\begin{proof}
 We refer to \cite[Theorem 2.4.63]{SE}.
\end{proof}

\begin{remark}
 Note that $\Delta$ in equation (\ref{BV-identity}) depends on which $m\in\mathbb{N}$ we choose to identify $\mathbb{D}:HH_\bullet(A)\stackrel{\sim}{\rightarrow} HH^{6m+5-\bullet}(A)[(2m+1)h+2]$, where the Gerstenhaber bracket does not.
\end{remark}

\subsection{Computation of the calculus structure of the preprojective algebra}
Since the calculus structure is defined on Hochschild chains and cochains, we have to work with the on the resolution for computations. It turns out that we only have to compute $\mathcal{L}_{\theta_0}$ directly, the rest can be deduced from formulas given by the calculus and the BV structure.

$$
\CD
\hdots@>d_3>>A\otimes A[2]@>d_2>>A\otimes V\otimes A@>d_1>>A\otimes A@>d_0>> A @>>> 0\\
@. @V\mu_2 VV @V\mu_1 VV @\vert @\vert\\
\hdots@>b_3>>A^{\otimes4}@>b_2>>A^{\otimes3}@>b_1>>A^{\otimes2}@>b_0>>A@>>>0
\endCD
$$
These maps $\psi_i$ give us a chain map between the Schofield and the bar resolution:
\begin{eqnarray*}
\mu_1(1\otimes y\otimes 1)&=&1\otimes y\otimes 1,\\
\mu_2(1\otimes 1)&=&\sum\limits_{a\in\bar Q}\epsilon_a 1\otimes a\otimes a^*\otimes1-1\otimes b\otimes b\otimes 1,\\
\mu_3(1\otimes1)&=&\sum\limits_{a\in\bar Q}\sum\limits_{x_i\in\mathfrak{B}}\epsilon_a 1\otimes \phi(x_i)\otimes a\otimes a^*\otimes x_i^*\\
&&-\sum\limits_{x_i\in\mathfrak{B}}1\otimes \phi(x_i)\otimes b\otimes b\otimes x_i^*,
\end{eqnarray*}
and 
\[
 \mu_{3+i}=\mu_i\left(\sum\limits_{a\in\bar Q}\sum\limits_{x_i\in\mathfrak{B}}\epsilon_a \phi(x_i)\otimes a\otimes a^*\otimes x_i^*-\sum\limits_{x_i\in\mathfrak{B}}\phi(x_i)\otimes b\otimes b\otimes x_i^*  \right).
\]

Now, we apply the functor $-\otimes_{A^e} A$ on the commutative diagram:

$$
\CD
\hdots@>d_3'>>A^R[2]@>d_2'>>(V\otimes A)^R @>d_1>>A^R@>>> 0 \\
@. @V\mu_2' VV @V\mu_1' VV @\vert \\
\hdots@>b_3'>>(A^{\otimes3})^R@>b_2>>(A^{\otimes2})^R@>b_1>>(A^{\otimes1})^R@>>>0
\endCD
$$

where

\begin{eqnarray*}
\mu_1'(x\otimes y)&=&x\otimes y,\\
\mu_2'(x)&=&\sum\limits_{a\in\bar Q}\epsilon_a a\otimes a^*\otimes x-b\otimes b\otimes x,\\
\mu_3'(x)&=&\sum\limits_{a\in\bar Q}\sum\limits_{x_i\in \mathfrak{B}}\epsilon_a \phi(x_i)\otimes a\otimes a^*\otimes x_i^*x-\phi(x_i)\otimes b\otimes b\otimes x_i^*x,
\end{eqnarray*}
and

\[
 \mu_{3+i}'=\mu_i'\left(\sum\limits_{a\in\bar Q}\sum\limits_{x_i\in \mathfrak{B}}\phi(x_i)\otimes a\otimes a^*\otimes x_i^*-\sum\limits_{x_i\in \mathfrak{B}}\phi(x_i)\otimes b\otimes b\otimes x_i^*\right).
\]

Now, we compute $\mathcal{L}_{\theta_0}$:

\begin{lemma}\label{Lie theta}
 For each $x\in HH_i(A)$, 
 \begin{equation}
  \mathcal{L}_{\theta_0}(x)=\deg(x)x
 \end{equation}
\end{lemma}

\begin{proof}
 Via $\mu'$, we already identified $x\in HH_i(A)$ with cycles in the Hochschild chain, but we still have to identify $\theta_0$ with an element in $\mathrm{Hom}_{A^e}(A^{\otimes3},A)$:

given any monomial $b=b_1\ldots b_l$, $b_i\in V$, the map 
\[
\tau(1\otimes b\otimes1)=\sum\limits_{i=1}^lb_1\ldots b_{i-1}\otimes b_i\otimes b_{i+1}\ldots b_l
\] 
makes the diagram
$$
\CD
A\otimes V\otimes A@>d_1>>A\otimes A@>d_0>> A @>>> 0\\
@A\tau AA @\vert @\vert\\
A^{\otimes3}@>b_1>>A^{\otimes2}@>b_0>>A@>>>0
\endCD
$$
commute.

Applying $\mathrm{Hom}_{A^e}(\_\_\otimes A)$, we get a map

\[
 \tau^*: \mathrm{Hom}_k(V)\rightarrow\mathrm{Hom}_k(A),
\]
such that
\[
 (\theta_0\circ\tau^*)(b_1\ldots b_l)=\sum\limits_{i=1}^lb_1\ldots b_{i-1}\theta_0(b_i) b_{i+1}\ldots b_l=\deg(b)\cdot b,
\]

Recall from \cite[(3.5), page 46]{D} that the Lie derivative of $\theta_0\circ\tau^*$ on Hochschild chains is defined by

\begin{eqnarray*}
 \mathcal{L}_{\theta_0\circ\tau^*}(a_1\otimes\cdots\otimes a_k)&=&\sum\limits_{i=1}^k a_1\otimes\cdots\otimes (\theta_0\circ\tau^*)(a_i)\otimes\cdots\otimes a_k\\
 &=&\sum\limits_{i=1}^k (\deg(a_1)+\cdots +\deg(a_k)) a_1\otimes\cdots\otimes a_k,
\end{eqnarray*}
and it can easily be checked that for each $x\in HH_i(A)$, $\mathcal{L}_{\theta_0\circ\tau^*}$ acts on $\mu_i'(x)$, $x\in HH^i(A)$, by multiplication with $\deg(x)$.

\end{proof}

\subsubsection{The contraction map}
From (\ref{intertwining}) we know that the contraction map on Hochschild homology is given by the cup product on Hochschild cohomology. Table \ref{contraction} contains these results, rewritten in terms of the contraction maps.

\subsubsection{The Connes differential}
We start with the computation of the Connes differential and refer the reader to the Subsection \ref{cyclic}.

\begin{proposition}\label{B}
 The Connes differential B is given by
\begin{eqnarray*}
 B_{6s}(\psi_{k,s})&=&((2s+1)h-2-k)\zeta_{k,s},\\
 B_{1+6s}&=&0,\\
 B_{2+6s}(h_{k,s})&=&(2s+1)h\alpha^{-1}(h_{k,s}),\\
 B_{3+6s}&=&0,\\
 B_{4+6s}(\theta_{k,s})&=&((2s+1)h+2+k)z_{k,s},\\
 B_{5+6s}&=&0.
\end{eqnarray*}
\end{proposition}
\begin{proof}
We use the Cartan identity (\ref{Cartan}) with $a\in\theta_0$,
\begin{equation}
\mathcal{L}_{\theta_0}=B\iota_{\theta_0}+\iota_{\theta_0} B,
\end{equation}
where $\mathcal{L}_{\theta_0}$ acts on $x\in HH_i$ by multiplication by $\deg(x)$ (see Lemma (\ref{Lie theta})).
The above identities for the Connes differential follow since $\iota_{\theta_0}$ acts on $\theta_{k,s}$, $\psi_{k,s}$ and $h_{k,s}$ by zero, and $z_{k,s}$,  $\zeta_{k,s}$ and $\alpha^{-1}(h_{k,s})$ are their unique preimages the contraction with $\iota_{\theta_0}$.
\end{proof}

\subsubsection{The Gerstenhaber bracket}

We compute the brackets using the identification \\$HH^{i}(A)=HH_{6m+5-i}(A)[-2(m+1)h-2]$ for $m>>1$ and the BV-identity (\ref{BV-identity}). We rewrite the results from Proposition \ref{B}:

\begin{eqnarray*}
 \Delta(\theta_k^{(s)})&=&((1+2(m-s))h+k+2)z_k^{(s)},\\
 \Delta(f_k^{(s)})&=&0,\\
 \Delta(h_k^{(s)})&=&(1+2(m-s))h\alpha^{-1}(h_k^{(s)}),\\
 \Delta(\zeta_k^{(s)})&=&0,\\
 \Delta(\psi_k^{(s)})&=&((1+2(m-s)h-k-2)\zeta_k^{(s)},\\
 \Delta(z_k^{(s)})&=&0.
\end{eqnarray*}

The cup products relations involving our basis of $HH^*(\Pi_{T_n})$ are the same ones as the relations in the $A_{2n}$-case. When comparing the differential $\Delta$ with the one in the $A_{2n}$-case where we identify $HH^{i}(\Pi_{A_{2n}})=HH_{6m+2-i}(\Pi_{A_{2n}})[-2mh-2]$ for $m>>1$, we have to multiply the coefficients by $2$ and add $h$. In the BV-identity (\ref{BV-identity}), we use only cup product and $\Delta$ to compute the Gerstenhaber bracket. In these computations, when comparing to the $A_{2n}$-case, we get the same results with the factor $2$. So using the results from \cite[Table 2]{Eu2}, we get Table \ref{Gerstenhaber bracket}.

\subsubsection{The Lie derivative $\mathcal{L}$}
We use the Cartan identity (\ref{Cartan}) to compute the Lie derivative.\\

\newpage

\textbf{\underline{$HH^{1+6s}(A)$-Lie derivatives}}:\\

From the Cartan identity, we see that
\[
\mathcal{L}_{\theta_k^{(s)}}=B\iota_{\theta_k^{(s)}}+\iota_{\theta_k^{(s)}}B.
\]
On $\theta_{l,t}$, $\psi_{l,t}$ and $h_{l,t}$, the Connes differential acts by multiplication with its degree and taking the preimage under $\iota_{\theta_0}$, and $\iota_{\theta_k^{(s)}}$ acts on them by zero. $B$ kills $z_{l,t}$, $\zeta_{l,t}$ and $f_{l,t}$. Since $B$ is degree preserving, this means that $\mathcal{L}_{\theta_k^{(s)}}$ acts on $\theta_{l,t}$, $\psi_{l,t}$ and $h_{l,t}$ by multiplication with their degree times $z_k^{(s)}$, and on $z_{l,t}$, $\zeta_{l,t}$ and $f_{l,t}$ by multiplication with $z_{k}^{(s)}$ and then multiplication with the degree of their product. So we get the following formulas:
\begin{eqnarray*}
 \mathcal{L}_{\theta_k^{(s)}}(\psi_{l,t})&=&((2t+1)h-2-l)(z_k\psi_l)_{t-s},\\ 
 \mathcal{L}_{\theta_k^{(s)}}(\zeta_{l,t})&=&((2(t-s)+1)h-2-l+k)(z_k\zeta_l)_{t-s},\\ 
 \mathcal{L}_{\theta_k^{(s)}}(h_{l,t})&=&\delta_{k0}(2t+1)hh_{l,t-s},\\
 \mathcal{L}_{\theta_k^{(s)}}(f_{l,t})&=&\delta_{k0}(2(t-s)+1)hf_{l,t-s},\\
 \mathcal{L}_{\theta_k^{(s)}}(\theta_{l,t})&=&((2t+1)h+2+l)(z_k\theta_l)_{t-s},\\
 \mathcal{L}_{\theta_k^{(s)}}(z_{l,t})&=&((2(t-s)+1)h+2+l+k)(z_kz_l)_{t-s}.
\end{eqnarray*}

\newpage

\textbf{\underline{$HH^{2+6s}(A)$-Lie derivatives}}:\\

We compute $\mathcal{L}_{f_k^{(s)}}$:
\begin{eqnarray*}
  \mathcal{L}_{f_k^{(s)}}(\psi_{l,t})&=&B(\iota_{f_k^{(s)}}(\psi_{l,t}))-\iota_{f_k^{(s)}}B(\psi_{l,t})\\
 &=&B(\delta_{l,h-3}k\theta_{h-3,t-s-1})-((2t+1)h-2-l)\iota_{f_k^{(s)}}(\zeta_{l,t})\\
 &=&\delta_{l,h-3}k((2(t-s-1)+2)h-1)z_{h-3,t-s-1}\\&&-\delta_{l,h-3}k((2t+1)h-2-l)z_{h-3,t-s-1}\\
 &=&-2\delta_{l,h-3}k((1+sh)z_{h-3,t-s-1},\\
  \mathcal{L}_{f_k^{(s)}}(\zeta_{l,t})&=&B(\iota_{f_k^{(s)}}(\zeta_{l,t}))=B(
k\delta_{l,h-3}z_{h-3,t-s-1})=0,\\
  \mathcal{L}_{f_k^{(s)}}(h_{l,t})&=&B(\iota_{f_k^{(s)}}(h_{l,t}))-\iota_{f_k^{(s)}}B(h_{l,t})\\
 &=&B(\delta_{k,l}\psi_{0,t-s})-(2t+1)h\iota_{f_k^{(s)}}\alpha^{-1}(h_{l,t})\\
 &=&\delta_{kl}((2(t-s)+1)h-2)\zeta_{0,t-s}-\delta_{kl}(2t+1)h\zeta_{0,t-s}\\
 &=&-2\delta_{kl}(sh+1)\zeta_{0,t-s},\\
 \mathcal{L}_{f_k^{(s)}}(f_{l,t})&=&B(\underbrace{\iota_{f_k^{(s)}}(f_{l,t})}_{\in HH_{1+6(t-s)}})=0,\\ 
 \mathcal{L}_{f_k^{(s)}}(\theta_{l,t})&=&B(\iota_{f_k^{(s)}}(\theta_{l,t}))-\iota_{f_k^{(s)}}(B(\theta_{l,t}))\\
 &=&B(\delta_{l0}\alpha(f_{k,t-s}))-((2t+1)h+2+l)\iota_{f_k^{(s)}}z_{l,t}\\
 &=&\delta_{l0}(2(t-s)+1)hf_{k,t-s}-\delta_{l0}((2t+1)h+2)f_{k,t-s}\\&=&-2\delta_{l0}(1+sh)f_{k,t-s},\\
 \mathcal{L}_{f_k^{(s)}}(z_{l,t})&=&\delta_{l0}B(f_{k,t-s})=0,
\end{eqnarray*}

\newpage

\textbf{\underline{$HH^{3+6s}(A)$-Lie derivatives}}:\\

We compute $\mathcal{L}_{h_k^{(s)}}$:
\begin{eqnarray*}
 \mathcal{L}_{h_k^{(s)}}(\psi_{l,t})&=&B(\underbrace{\iota_{h_k^{(s)}}(\psi_{l,t})}_{=0})+\iota_{h_{k}^{(s)}}B(\psi_{l,t})\\
 &=&((2t+1)h-2-l)\iota_{h_k^{(s)}}\zeta_{l,t}=\delta_{k,n}\delta_{l,h-3}(2th+1)\theta_{h-3,t-s-1},\\
 \mathcal{L}_{h_k^{(s)}}(\zeta_{l,t})&=&B(\iota_{h_k^{(s)}}(\zeta_{l,t}))=\delta_{k,n}B(\delta_{l,h-3}\theta_{h-3,t-s-1})\\&=&\delta_{kn}\delta_{l,h-3}((2(t-s-1)+1)h+2+h-3)z_{h-3,t-s-1}\\
 &=&\delta_{kn}\delta_{l,h-3}((2(t-s)h-1)z_{h-3,t-s-1},\\
 \mathcal{L}_{h_k^{(s)}}(h_{l,t})&=&B(\underbrace{\iota_{h_k^{(s)}}(h_{l,t})}_{=0})+\iota_{h_k^{(s)}}B(h_{l,t})=(2t+1)h\iota_{h_k^{(s)}}\alpha^{-1}(h_{l,t})\\
 &=&(2t+1)h(M_\alpha^{-1})_{lk}\psi_{0,t-s}\\
 \mathcal{L}_{h_k^{(s)}}(f_{l,t})&=&B(\iota_{h_k^{(s)}}(f_{l,t}))=B(\delta_{kl}\psi_{0,t-s})=\delta_{kl}(2(t-s+1)h-2)\zeta_{0,t-s},\\ 
 \mathcal{L}_{h_k^{(s)}}(\theta_{l,t})&=&B(\underbrace{\iota_{h_k^{(s)}}(\theta_l^{(t)})}_{=0})+\iota_{h_k^{(s)}}B(\theta_{l,t})=((2t+1)h+2+l)\iota_{h_k^{(s)}}z_{l,t}\\
 &=&\delta_{l0}((2t+1)h+2)h_{k,t-s},\\
 \mathcal{L}_{h_k^{(s)}}(z_{l,t})&=&B(\delta_{l0}h_{k,t-s})=\delta_{l0}(2(t-s)+1)h\alpha^{-1}(h_{k,t-s}).
\end{eqnarray*}

\newpage

\textbf{\underline{$HH^{4+6s}(A)$-Lie derivatives}}:\\

We compute $\mathcal{L}_{\zeta_k^{(s)}}$:
\begin{eqnarray*}
 \mathcal{L}_{\zeta_k^{(s)}}(\psi_{l,t})&=&B\iota_{\zeta_k^{(s)}}(\psi_{l,t})-\iota_{\zeta_k^{(s)}}B(\psi_{l,t})\\
 &=&\delta_{k,h-3}\delta_{l,h-3}B(\alpha(f_{n,t-s-1}))-\iota_{\zeta_k^{(s)}}((2t+1)h-2-l)\zeta_{l,t}\\
 &=&\delta_{k,h-3}\delta_{l,h-3}((2(t-s-1)+1)hf_{n,t-s-1}\\&&-((2t+1)h-h+1)f_{n,t-s-1})\\
 &=&\delta_{k,h-3}\delta_{l,h-3}((-2s-1)h-1)f_{n,t-s-1},\\
 \mathcal{L}_{\zeta_k^{(s)}}(\zeta_{l,t})&=&B\iota_{\zeta_k^{(s)}}(\zeta_{l,t})-\iota_{\zeta_k^{(s)}}\underbrace{B(\zeta_{l,t})}_{=0}\\
 &=&\delta_{k,h-3}\delta_{l,h-3}B(f_{n,t-s-1})=0,\\
 \mathcal{L}_{\zeta_k^{(s)}}(h_{l,t})&=&B\iota_{\zeta_k^{(s)}}(h_{l,t})-\iota_{\zeta_k^{(s)}}B(h_{l,t})\\
 &=&
 \delta_{l,n}\delta_{k,h-3}B(\theta_{h-3,t-s-1})-(2t+1)h\iota_{\zeta_k^{(s)}}\alpha^{-1}(h_{l,t}),\\
 &=&\delta_{l,n}\delta_{k,h-3}z_{h-3,t-s-1}((2(t-s-1)+1)h+2+h-3-(2t+1)h)\\
 &=&\delta_{l,n}\delta_{k,h-3}z_{h-3,t-s-1}(-(2s+1)h-1),\\
 \mathcal{L}_{\zeta_k^{(s)}}(f_{l,t})&=&B\iota_{\zeta_k^{(s)}}(f_{l,t})-\iota_{\zeta_k^{(s)}}\underbrace{B(f_{l,t})}_{=0}\\
 &=&l\delta_{k,h-3}B(z_{h-3,t-s})=0,\\
 \mathcal{L}_{\zeta_k^{(s)}}(\theta_{l,t})&=&B\iota_{\zeta_k^{(s)}}(\theta_{l,t})-\iota_{\zeta_k^{(s)}}B(\theta_{l,t})=B((z_l\psi_k)_{t-s})-\iota_{\zeta_k^{(s)}}((2t+1)h+2+l)z_{l,t}\\
 &=&((2(t-s)+1)h-2-(k-l))(z_l\zeta_k)_{t-s}-((2t+1)h+2+l)(z_l\zeta_k)_{t-s-1}\\
 &=&(-2sh-4-k)(z_l\zeta_k)_{t-s},\\
 \mathcal{L}_{\zeta_k^{(s)}}(z_{l,t})&=&B\iota_{\zeta_k^{(s)}}(z_{l,t})-\iota_{\zeta_k^{(s)}}\underbrace{B(z_{l,t})}_{=0}\\
 &=&B((z_l\zeta_k)_{t-s})=0.
\end{eqnarray*}

\newpage 

\textbf{\underline{$HH^{5+6s}(A)$-Lie derivatives}}:\\

We compute $\mathcal{L}_{\psi_k^{(s)}}$:

\begin{eqnarray*}
 \mathcal{L}_{\psi_k^{(s)}}(\psi_{l,t})&=&B\underbrace{\iota_{\psi_k^{(s)}}(\psi_{l,t})}_{=0}+\iota_{\psi_k^{(s)}}B(\psi_{l,t})\\ 
 &=&((2t+1)h-2-l)\iota_{\psi_k^{(s)}}\zeta_{l,t}\\
 &=&\delta_{k,h-3}\delta_{l,h-3}(\underbrace{(2t+1)h-2-(h-3)}_{=2th+1}\alpha(f_{n,t-s-1}),\\
 \mathcal{L}_{\psi_k^{(s)}}(\zeta_{l,t})&=&B\iota_{\psi_k^{(s)}}(\zeta_{l,t})+\iota_{\psi_k^{(s)}}\underbrace{B(\zeta_{l,t})}_{=0}\\
 &=&\delta_{k,h-3}\delta_{l,h-3}B(\alpha(f_{n,t-s-1}))\\
 &=&\delta_{k,h-3}\delta_{l,h-3}(2(t-s-1)+1)hf_{n,t-s-1}\\
 &=&\delta_{k,h-3}\delta_{l,h-3}(2(t-s)-1)hf_{n,t-s-1},\\
 \mathcal{L}_{\psi_k^{(s)}}(h_{l,t})&=&B\underbrace{\iota_{\psi_k^{(s)}}(h_{l,t})}_{=0}+\iota_{\psi_k^{(s)}}B(h_{l,t})\\
 &=&\iota_{\psi_k^{(s)}}\alpha^{-1}(h_{l,t})(2t+1)h\\
 &=&\delta_{k,h-3}\delta_{l,n}(2t+1)h\theta_{h-3,t-s-1},\\
 \mathcal{L}_{\psi_k^{(s)}}(f_{l,t})&=&B\iota_{\psi_k^{(s)}}(f_{l,t})+\iota_{\psi_k^{(s)}}\underbrace{B(f_{l,t})}_{=0}\\
 &=&l\delta_{k,h-3}B(\theta_{h-3,t-s-1})\\&=&l((2(t-s-1)+1)h+2+(h-3))\delta_{k,h-3}z_{h-3,t-s}\\
 &=&l(2(t-s)h-1)\delta_{k,h-3}z_{h-3,t-s-1},\\
 \mathcal{L}_{\psi_k^{(s)}}(\theta_{l,t})&=&B\underbrace{\iota_{\psi_k^{(s)}}(\theta_{l,t})}_{=0}+\iota_{\psi_k^{(s)}}B(\theta_{l,t})\\
 &=&\iota_{\psi_k^{(s)}}z_{l,t}((2t+1)h+2+l)=(z_l\psi_k)_{t-s}((2t+1)h+2+l),\\
 \mathcal{L}_{\psi_k^{(s)}}(z_{l,t})&=&B\iota_{\psi_k^{(s)}}(z_{l,t})+\iota_{\psi_k^{(s)}}\underbrace{B(z_{l,t})}_{=0}\\
 &=&B((z_l\psi_k)_{t-s})=((2(t-s)+1)h-2-(k-l))(z_l\zeta_k)_{t-s}.\\
\end{eqnarray*}

\newpage

\textbf{\underline{$HH^{6+6s}(A)$-Lie derivatives}}:\\

$B$ acts on $\theta_{l,t},\psi_{l,t}$ and $h_{l,t}$ by multiplication with its degree and taking the preimage under $\iota_{\theta_0}$. On $z_{l,t},\zeta_{l,t}$ and $f_{l,t}$, $B$ acts by zero. Since the spaces $U$, $U^*$, $K$ and $K^*$ are $z_k$-invariant and $z_k^{(s)}$ has degree $k-2sh$, $\mathcal{L}_{z_k^{(s)}}$ acts on  $\theta_{l,t},\psi_{l,t}$ and $h_{l,t}$ by multiplication with $k-2sh$ and taking the preimage under $\iota_{\theta_0}$ and multiplication with $z_k^{(s)}$, and on $z_{l,t},\zeta_{l,t}$ and $f_{l,t}$ it acts by zero. We have the following formulas:

\begin{eqnarray*}
 \mathcal{L}_{z_k^{(s)}}(\psi_{l,t})&=&(k-2sh)(z_k\zeta_l)_{t-s},\\
 \mathcal{L}_{z_k^{(s)}}(\zeta_{l,t})&=&0,\\
 \mathcal{L}_{z_k^{(s)}}(h_{l,t})&=&(k-2sh)\alpha^{-1}(h_{l,t-s}),\\
 \mathcal{L}_{z_k^{(s)}}(f_{l,t})&=&0,\\
 \mathcal{L}_{z_k^{(s)}}(\theta_{l,t})&=&(k-2sh)(z_kz_l)_{t-s},\\
 \mathcal{L}_{z_k^{(s)}}(z_{l,t})&=&0.\\
\end{eqnarray*}

\appendix
\section{Correction to \cite{ES2}}
We want to make a correction to the $HH^2(A)\cup HH^2(A)$-computation in \cite{ES2}: the calculation of $HH^2(A)\cup HH^2(A)$ in \cite[Subsection 16.2]{Eu2} shows that the bilinear form on $K$ is given by the matrix $M_\alpha$, defined in \cite[Subsection 15.2.]{Eu2}. This is a general computation which also applies to quivers of type A. But the results in \cite{ES2} suggest that the bilinear form on $K$ is given by a matrix different from $M_\alpha$ which is incorrect. 

I verified that the matrix $M_\alpha$ from $HH^1(A)\cup HH^2(A)$ in \cite{ES2} correct, therefore the result of $HH^2(A)\cup HH^2(A)$ is wrong: similarly to the computation in Subsection \ref{H1xH2} of this paper, you can calculate the matrix $M_\alpha$ by using the derivative of $H_A(t)$. Then you get (by labeling the $A$-quiver as in \cite{ES2})

\[
 M_\alpha=
 h
 \left[
 \begin{array}{cccccc}
 2&-1&0&\cdots&\cdots&0\\
 -1&2&\ddots&\ddots&&\vdots\\
 0&\ddots&\ddots&\ddots&\ddots&\vdots\\
 \vdots&\ddots&\ddots&\ddots&\ddots&0\\
 0&\cdots&0&-1&2&-1\\
 0&\cdots&\cdots&0&-1&2
 \end{array}
 \right]^{-1}\]

for type $A_{2n+1}$ (and also $D_n,\, E_n$)
and 
\[
 M_\alpha=
 h
 \left[
 \begin{array}{cccccc}
 2&-1&0&\cdots&\cdots&0\\
 -1&2&\ddots&\ddots&&\vdots\\
 0&\ddots&\ddots&\ddots&\ddots&\vdots\\
 \vdots&\ddots&\ddots&\ddots&\ddots&0\\
 0&\cdots&0&-1&2&-1\\
 0&\cdots&\cdots&0&-1&3
 \end{array}
 \right]^{-1}\]
for type $A_{2n}$ which are exactly the same results as in \cite{ES2}.

\end{document}